\documentclass{article}
\usepackage{graphicx} 
\graphicspath{{./}}
\usepackage{subcaption}
\usepackage{amsmath}
\usepackage{xcolor}
\usepackage{biblatex}
\usepackage{authblk}
\usepackage[section]{placeins} 
\usepackage{multirow}
\usepackage{verbatim}
\title{GO Competition Challenge 3: Problem, Solvers, and Solution Analysis}
\author[1]{Jesse T. Holzer\thanks{jesse.holzer@pnnl.gov}}
\author[1]{Stephen Elbert\thanks{steve.elbert@pnnl.gov}}
\author[2]{Hans Mittelmann\thanks{mittelma@asu.edu}}
\author[3]{Richard O'Neill\thanks{richardponeill@gmail.com}} 
\author[4]{HyungSeon Oh\thanks{hyungseon.oh@hq.doe.gov, corresponding author}}
\affil[1]{Pacific Northwest National Laboratory}
\affil[2]{Arizona State University}
\affil[3]{Advanced Research Projects Agency - Energy}
\affil[4]{Booz Allen Hamilton}
\date{August 2024}

\begin{document}

\maketitle

\begin{abstract}
    This paper describes the Grid Optimization (GO) Competition Challenge 3, focusing on the problem motivation, formulation, solvers submitted by competition entrants, and analysis of the solutions produced. Funded by DOE/ARPA-E and led by a collaboration of national labs and academia members, the GO Competition addresses challenging problems in power systems planning and operations to drive research in advanced solution methods essential for a rapidly evolving electric power sector. Challenge 3 targets a multi-period unit commitment problem, incorporating AC power modeling and topology switching to reflect the dynamic grid management techniques required for future power systems. The competition results offer significant benefits to both researchers and industry practitioners. For researchers, it fosters innovation, encouraging the development of new algorithms to address the complexities of modern power systems. For industry practitioners, the competition drives the creation of more efficient and reliable computational tools, directly improving grid management practices. This collaboration bridges the gap between theory and practical implementation, advancing the field in meaningful ways. This paper documents the problem formulation, solver approaches, and the effectiveness of the solutions developed.
\end{abstract}

\section*{Acknowledgments}

This paper would not have been possible without the major contributions of the entire Grid Optimization (GO) Competition support team, both directly in the development of the paper and indirectly in the work of running the competition. Over the course of the GO Competition project, the support team has included the following individuals, in alphabetical order:
  Carleton Coffrin,
  Christopher DeMarco,
  Ray Duthu,
  Stephen Elbert,
  Brent Eldridge,
  Tarek Elgindy,
  Manuel Garcia,
  Scott Greene,
  Nongchao Guo,
  Elaine Hale,
  Jesse Holzer,
  Bernard Lesieutre,
  Terrence Mak,
  Colin McMillan,
  Hans Mittelmann,
  HyungSeon Oh,
  Richard O'Neill,
  Thomas Overbye,
  Bryan Palmintier,
  Robert Parker,
  Farnaz Safdarian,
  Ahmad Tbaileh,
  Pascal Van Hentenryck,
  Arun Veeramany,
  Steve Wangen,
  Jessica Wert.
 The GO Competition project has been supported from its inception by the United States Department of Energy (DOE) Advanced Research Projects Agency - Energy (ARPA-E) award number 68150 "Supporting ARPA-E Power Grid Optimization".

\section{Introduction}

This paper describes the solvers developed by entrants in the Grid Optimization (GO) Competition Challenge 3.
The Grid Optimization Competition is run by a team of researchers from a number of organizations, including the sponsor Advanced Research Projects Agency - Energy (ARPA-E), lead organization Pacific Northwest National Laboratory (PNNL), and technical contributors from Los Alamos National Laboratory (LANL), National Renewable Energy Laboratory (NREL), Texas A\&M University (TAMU), Georgia Institute of Technology (GT), University of Wisconsin (UW), and others.
The GO Competition poses challenge problems in the field of power grid management, invites entrants to develop solvers for these problems, invokes the solvers on a set of problem instances using common hardware, ranks the solvers according to their performance, and awards prizes according to the rankings.
The overall goal of the GO Competition is to spur innovative research on high impact and computationally challenging problems in power grid management from initial development through commercial deployment.

The GO Competition has run three challenges. This paper covers Challenge 3. The main objective of the paper is to describe the solution methods developed by the entrants and the competition organizers and to analyze the performance of these different solvers. 
In Section \ref{sec:prob_significance}, we describe the significance of the Challenge 3 problem and its technical difficulty.
In Section \ref{sec:prob_desc}, we describe the Challenge 3 problem in greater detail, including our motivation for certain features. In Section \ref{sec:lit_rev}, we review the literature relevant to this problem. In Section \ref{sec:prob_form}, we give a condensed formulation of the problem. In Section \ref{sec:sol_methods_traditional}, we review traditional solution approaches. In Section \ref{sec:sol_methods_entrants}, we describe the approaches of the competition entrants. In Section \ref{sec:results}, we present some results of the competition. In Section \ref{sec:conclusion}, we offer conclusions.

\section{Importance, Challenges, and Implications of Solving the GO Competition Challenge 3 Problem of Multi-Period Unit Commitment with AC Power Flow}
\label{sec:prob_significance}

Combining unit commitment (UC) with a full alternating current optimal power flow (AC OPF) model is a highly sophisticated and complex challenge in power system operations. However, the benefits of solving this problem are significant, offering enhanced accuracy and reliability in managing the electrical grid.

The primary advantage of incorporating a full AC OPF model into unit commitment lies in the increased accuracy of the AC power flow model. Unlike the more commonly used direct current (DC) power flow model, which assumes constant voltage magnitudes and neglects reactive power, the full AC model represents the power flow equations in true nonlinear form. This allows for precise modeling of voltage magnitudes, angles, losses, and reactive power, which are crucial for making accurate operational decisions, particularly in systems where voltage stability is a key concern. The nonlinearity of the AC power flow equations also allows a more realistic representation of power transmission, including losses that are often ignored or approximated in DC models. The linear approximations made by a DC power flow model imply that, in order to use these models, the modeler must either accept that the solution will sometimes violate certain engineering constraints such as line flow limits and voltage constraints or must enforce these constraints by using conservative limits sometimes leading to suboptimal solutions. Using an AC power flow model avoids this tradeoff between violating engineering constraints and suboptimal dispatch, enabling bette economic performance without compromising power system security.

Moreover, the reliability and stability of the power grid are greatly enhanced when using a full AC OPF model. Voltage stability is more accurately managed, which is critical for maintaining the grid's reliability, particularly during times of high demand or stress. The improved accuracy in modeling leads to better operational decisions, to avoid costly or potentially unstable conditions that might arise from the approximations inherent in DC OPF models. This is particularly important as power systems increasingly integrate renewable energy sources, which are variable and can have significant impacts on both active and reactive power flows.

However, the complexity of solving UC combined with a full AC OPF model presents several significant challenges. One of the main difficulties is the inherent nonlinearity of the AC power flow equations. These equations are nonlinear and nonconvex, which makes the problem much more complex than the linear DC OPF model. This nonlinearity leads to the presence of multiple local minima, making it difficult to find a globally optimal solution. When unit commitment, which involves making discrete decisions such as turning generators on or off, is combined with the nonlinear AC OPF model, the result is a mixed-integer nonlinear programming (MINLP) problem. This is vastly more challenging to solve than the mixed-integer linear programming (MILP) problems typically used with DC OPF models. As the size of the power system increases, the computational burden grows exponentially, making it difficult to solve these problems in real-time or even near real-time, which is necessary for day-ahead or intra-day market operations.

Another challenge is the issue of convergence. Solving the AC OPF component typically requires iterative methods, such as Newton-Raphson or interior-point methods, which can struggle with convergence, particularly in systems that are heavily loaded or under stress.  This issue is further complicated by the binary decision-making process involved in unit commitment, making it even harder to find a stable and optimal solution.

Current software and algorithmic tools are often optimized for either MILP (used in UC with DC OPF) or nonlinear continuous optimization (used in AC OPF alone), but not for the combined MINLP problem that arises when both are integrated. This gap in existing tools underscores the need to develop or to improve algorithms that can efficiently handle this combined problem.

Despite these challenges, the significance of overcoming them is profound. Successfully solving the UC problem with a full AC OPF model would lead to enhanced grid reliability by ensuring more accurate and stable operation, thus reducing the risk of blackouts or other reliability issues related to voltage instability. Economically, more accurate modeling could lead to more efficient use of generation resources, potentially lowering the overall cost of electricity supply. This approach would also facilitate the integration of higher levels of renewable energy, which is critical as the world transitions to sustainable energy sources. As electricity markets and regulations evolve, particularly with stricter emissions and renewable energy mandates, the ability to solve this combined problem becomes increasingly critical for market participants to remain compliant and competitive.

In conclusion, while integrating unit commitment with a full AC OPF model presents significant technical and computational challenges, overcoming these hurdles would represent a major advancement in power system operations. This would enable more accurate, reliable, and economically efficient management of the electrical grid, particularly in the context of growing renewable energy integration and evolving market dynamics.

\section{Problem Description}
\label{sec:prob_desc}

In this section, we describe the GO Competition Challenge 3 problem at a high level and explain our motivation for including certain features.


In an uncertain era, the evolution of power systems through the integration of renewable energy sources and the deployment of advanced smart grid technologies poses a complex challenge. The GO Competition Challenge 3 problem is at the heart of this transformation, focusing on the optimization of short-term planning in electrical power systems that must now deal with the unpredictable and uncontrollable nature of renewable energy. This challenge highlights an important aspect of modern power systems: as we get closer to execution time, the accuracy of resource availability and weather forecasts improves dramatically, providing a critical opportunity to refine unit commitment strategies in real time.

Our research focuses on the efforts to integrate renewable energy into the grid, recognizing its environmental benefits despite its intermittent and non-controllable nature. By focusing on the transmission level and involving a large network of suppliers and consumers, we advocate for a comprehensive management approach that not only addresses immediate operational needs but also anticipates the changing landscape of energy generation and consumption.

The Grid Optimization challenge we investigate is intended to maximize efficiency and market surplus while taking into account the stochastic nature of renewable resources and the increasing accuracy of forecasts. This problem, traditionally referred to as an alternating current (AC) security-constrained (SC) unit commitment (UC) problem, becomes more complex as renewables are integrated. However, advances in forecasting technologies provide a distinct advantage; as the execution date approaches, improved precision in predicting resource availability and weather conditions can be seamlessly integrated into unit commitment decisions. This adaptive approach enables more accurate planning by optimizing generator startup/shutdown and power dispatch processes to account for the real-time characteristics of renewable energy sources.

Our contribution on this subject is focused on innovative analytics and methodologies that exploit the improving accuracy of forecasts to navigate the integration challenges of renewable energy and smart grid technologies. Highlighting strategies that dynamically adjust to the enhanced forecast precision by real-time UC decisions, our work aims to showcase how power systems can achieve resilience, competence, and sustainability. By embracing the temporal improvement in data accuracy, we propose solutions that not only ensure the stability and reliability of the grid but also secure the way for a more sustainable energy future.


The GO Competition Challenge 3 problem is an optimization problem for short term planning of an electrical power system with special model features intended to reflect the future needs of grid planners in a rapidly changing electric power sector. We consider the power system at the transmission level, with suppliers such as generators, and consumers such as local utilities and large industrial facilities, all connected by a network of transmission lines. In the overall application of short term planning, we are focused on decisions that might be made several minutes to several days ahead of time, mainly which generators to start up or shut down how much power to dispatch from each one, but we also consider decisions such as load dispatch and line switching that might become more salient in a future power grid. By casting the problem as an optimization problem,
we are looking for a plan that is the best among all possible plans
according to a defined objective such as the total cost of generation. More generally, we consider dispatchable load, so we also include the total value of load served, and therefore we optimize the total market surplus.

This problem would traditionally be called an alternating current (AC) security constrained (SC) unit commitment (UC) problem. UC means we are considering the commitment of generating units through discrete decisions on startup and shutdown as well as their dispatch through continuous decisions on real power output. AC modeling covers the dispatch, flow, and balance of reactive power in addition to real power, along with voltage magnitude. SCs are constraints that ensure the planned operating point stays within safe operating limits in normal conditions and credible outage contingencies.


We consider a UC problem in a typical multi-period context, with discrete time intervals covering a planning horizon of a few hours up to ten or more days. We specifically identify three applications based on the time horizon.
A 24 to 48 hour horizon with 15 minute to 2 hour time
intervals is used for the day ahead planning context common in wholesale electricity markets.
A 4 to 8 hour horizon with 15 minute to 1 hour time intervals is used for the near real time look ahead context where, for example, additional resources are brought online to handle revised forecasts of wind and solar output.
A 5 to 10 day horizon with 1 hour to 6 hour intervals is used for a week ahead planning process that might be used especially to prepare for severe weather events. This approach helps create a reliable operation plan to respond stochastic weather conditions like abnormal low temperatures in Texas in 2021 and the Eastern United States in 2023.

In UC, the discrete decisions of generator startup and shutdown are modeled with binary variables, so UC is a problem of mixed integer programming (MIP). As a MIP problem UC is NP hard. This implies that solution algorithms can take too long to reach a solution with a proof of optimality within a prescribed tolerance. The practical performance of modern MIP solvers on UC is typically much better than the theoretical worst case, but there is no guarantee of solver performance, and long run times can and sometimes do happen. This poses an ongoing challenge for UC applications.

AC means we consider not only the real power output of generators and the physics of real power balance and flow but also reactive power and voltage constraints. AC modeling is typically not used in short term planning in current practice, but there is reason to believe it can permit more efficient use of generation and transmission resources in the current grid and that it may become more important as the power grid evolves.

In general, after commitment and dispatch decisions are made, the state of the power grid is realized and follows physical laws of AC power flow and balance. These physical laws are formulated as nonlinear equations in the variables of device-level real and reactive power and bus-level voltage magnitude and angle. As an optimization problem, the optimal dispatch of generators and loads (even fixed loads) subject to AC physics is a nonconvex nonlinear programming (NLP) problem. In terms of computational difficulty this problem is NP-hard.

Without the physical laws of AC power flow and balance, it is not possible to model the bus voltage magnitudes and reactive power capabilities and requirements of devices. In particular the DC power flow and balance model typically used in UC and economic dispatch models in power grid operational planning at the time scales we consider here cannot model voltage and reactive power. In practice, when a dispatch solution obtained from a DC model is implemented, the resulting bus voltages may violate engineering limits, and real time control mechanisms keep the voltages in acceptable ranges. These control mechanisms depend on the reactive power and other capabilities and requirements of devices such as generators, shunts, and loads. If these capabilities are used to their limits, then the real time controls will no longer be able to maintain voltages within the desired ranges, and damage to grid equipment and cascading grid failure can result. If reactive power and voltage were considered at the dispatch planning stage by incorporating AC power flow and balance in the UC and economic dispatch models, then the likelihood of this failure mode could be decreased. In particular, additional generators might be committed in order to make their reactive power capability available. Furthermore, with reactive power capability characteristics (for example D-curves), the real power dispatch might be modified in order to enable any given device to provide more reactive power.

In the past, it has been unnecessary to consider voltage and reactive power at the dispatch planning stage because it has been possible to ensure sufficient reactive power capability mostly by relying on static reactive power infrastructure such as shunts but also by applying special constraints to ensure certain generators are committed in thoroughly studied and well understood situations of voltage stress and by generally applying conservative flow limits to transmission lines. We believe that changes underway in the electric power industry will pose a challenge to this method of handling voltage and reactive power after planning the real power dispatch. Specifically, the geographical and temporal variability of wind and solar power and the load flexibility that might be needed to manage this variability will lead to a much greater diversity of dispatch conditions and operating points, so that fixed infrastructure and deep offline study of all credible voltage stress scenarios will no longer be practical. Furthermore,
even in the current power grid,
AC modeling will permit the use of less conservative flow limits on transmission lines allowing in turn more efficient dispatch and commitment decisions.

Therefore, we include a full AC power flow and balance model at each time step of the UC model. In terms of theoretical computational complexity, UC and ACOPF are already hard problems on their own, and combining them does not change this theoretical complexity. However, UC and ACOPF on their own can each be modeled and solved with solver software at a high level of commercial and academic maturity. Commercial solvers for MIP have benefited from decades of intense development spurred by high value applications throughout modern society, and mature and robust NLP solvers are also available. Practical instances of UC and ACOPF can be solved in this way with reasonable performance. Solvers for the combined MINLP problem are at an early stage of development, so that modeling a combined UC-AC problem in a straightforward fashion and passing the model to a general purpose MINLP solver is not successful on any but the smallest problem instances. With this increased practical difficulty, in addition to the theoretical complexity, this problem fits the ARPA-E model of spurring research on the hardest problems.

SC refers to a wide variety of constraints ensuring not only that the dispatch plan is consistent with physical laws of electricity but also that the resulting system state remains within safe operating limits established by engineering practice. In the UC context, SCs can be viewed narrowly as the constraints ensuring that power flows on transmission lines do not exceed predetermined limits either in the base case when the network equipment is all in operation or in any of a set of credible contingencies each defined by the unplanned outage of one or more pieces of transmission equipment.

More generally, the concept of security to credible contingencies also requires that we consider the unplanned outage of generation equipment, and our model handles this in the traditional way by requiring reserves of generation capacity that could be called up in case of a generator outage to replace the power that was being provided by that outaged generator. Furthermore, the concept of safe operation requires that we consider voltage limits with AC modeling. Finally, for simplicity, our post-contingency model for line outage contingencies is a real power only model and thus cannot resolve post-contingency bus voltages, so we introduce reactive power reserve requirements to ensure that safe voltages can be maintained in a contingency.

The future power grid is expected to be more reliant on wind and solar, which are attractive for their near zero incremental energy cost and greenhouse gas emissions but that also have high variability in available power output. Therefore, it is anticipated that dispatchable load will be critical to maintaining the balance of energy supply and demand at all times. Furthermore, some of the types of loads that are likely to have the most potential for dispatchability are large industrial plants, such as metallurgical smelters, cement manufacturers, chemical processors and refiners, and even carbon capture and sequestration equipment.
Such loads may have significant operational complexity analogous to the startup, shutdown, and minimum uptime and downtime constraints of generators that are familiar in UC.

We therefore model dispatchable loads with all the same features as generators. From the standpoint of computational complexity and practical algorithmic performance, the main dimensions determining the computational difficulty of a UC problem are the number of generators and the number of time intervals. Therefore, this model feature reflecting a future user need might transform a fairly difficult but reasonable UC problem of 1000 generators and 24 or 48 time intervals into an enormously difficult problem of 5000 generator-like producing and consuming devices.

A further implication of the increasing reliance on wind and solar is that the topology of the power grid may need to be changed frequently in response to weather conditions in order to take best advantage of the available wind and solar energy. The greater diversity of dispatch conditions caused by the variability of wind and solar means that some lines should be switched (i.e. either connected or disconnected) in order to permit more power flow overall. We therefore include a decision variable to open or close each line in the network at each time step. In current practice, lines are typically not opened or closed due to day ahead or week ahead planning, and intense offline study is required to ensure that such topology switching actions can be performed without adversely affecting the dynamic state of the grid. We believe that much of the line switching analysis could be brought into the daily and weekly planning processes. Our competition takes a step towards doing that by including topology switching in the formulation to investigate the value that topology optimization could provide.

\section{Literature Survey}
\label{sec:lit_rev}

{The GO Challenge 3 delves into the complexities of power systems, where Unit Commitment (UC) and AC Optimal Power Flow (OPF) represent critical yet distinct aspects of modern electricity management. While UC primarily uses convexified (typically linear) power flow equations to schedule generation units, AC OPF typically operates with predetermined generation commitments. These topics are tailored to their unique characteristics and constraints, recognizing that the territory we explore has not been extensively charted in previous research.

The evolution of UC has been marked by a significant shift from deterministic to dynamic models, which corresponds to the complexities of evolving electricity markets and grid infrastructures. In its early stages, UC relied heavily on traditional methods such as Priority-List, Dynamic Programming, and Lagrangian Relaxation to optimize the scheduling of generation units. These methods, which were fundamental in their inception, were investigated in early studies \cite{4501405}. 
As renewable energy sources became more prevalent in the energy sector, UC was forced to adapt, resulting in the development of stochastic and robust models capable of dealing with the new layers of variability and uncertainty \cite{4470561}. 
This time period also saw the introduction of advanced optimization techniques, particularly Mixed-Integer Linear Programming (MILP), which signaled the beginning of a new era in UC research focused on improving solution accuracy and efficiency \cite{5983423}. 
Researchers developed decomposition and hybrid methods, such as Benders Decomposition \cite{1425598}, 
to deal with the growing scale of UC problems. Recent emphasis has been placed on integrating renewable energy, demand response, and storage, indicating a shift toward more sustainable power systems \cite{867149} 
Theoretical advances in UC have been matched by advances in specialized software and tools, such as MATPOWER, which have significantly improved computational efficiency and practical application \cite{5491276}. 
According to recent IEEE publications \cite{O’Neill2001}, 
the current trajectory of UC research is firmly oriented towards real-time algorithm development.

The journey in the field of ACOPF began with a focus on formulating and solving nonlinear optimization problems to minimize generation costs while adhering to operational constraints, laying the groundwork for future research \cite{carpentier1962}. 
Researchers shifted their focus to developing more efficient and robust solutions to address the nonlinearity, nonconvexity, and complexity inherent in power system models as computing power improved \cite{4075418}. 
The incorporation of uncertainty and variability, largely driven by the growing influence of renewable energy sources, was a critical moment in ACOPF research. As emphasized in \cite{CAPITANESCU20111731}, 
robust optimization techniques were required to ensure resilient power system operations. Power market deregulation and restructuring further complicated the ACOPF landscape, bringing market dynamics to the forefront of OPF formulations and solutions \cite{bialek1996}. 
Recent studies have incorporated renewable energy and demand response into ACOPF models to keep up with the evolving power grid \cite{conejo2010}. 
Along with theoretical advances, computational advances have significantly improved the scalability and efficiency of ACOPF solutions \cite{7271127}. 
The domain is now on the verge of a new era, with machine learning and artificial intelligence increasingly being used to augment the predictive and optimization capabilities of ACOPF models, indicating a promising future \cite{9042547}. 

As a result, our investigation in this study contributes these two critical areas of power systems - UC and ACOPF - with a recognition of their distinct yet interconnected nature. We hope to deepen a more nuanced understanding of modern power systems by addressing these issues separately but within a unified framework, guided by historical context and driven by advancements and future potential in these fields.

}

\section{Problem Formulation}
\label{sec:prob_form}
This section gives an abbreviated formulation of the model solved by the GO Competition. The fully detailed formulation used by the GO Competition, which we refer to here as the full formulation, is given in \cite{go_competition_ch_3_formulation}. Here we use complex numbers for compact expression of AC power flow. $\sqrt{-1}$ is the imaginary unit, and $z^*$ denotes the complex conjugate of $z$. We also use the indicator function notation $1_Z$ denoting the value $1$ if $Z$ is true and $0$ else.
%
%
\subsection{Nomenclature}
%
Each symbol used in the formulation is described as it is introduced, but we provide tables of nomenclature in this section for the convenience of the reader.

\subsubsection{Index sets}
\begin{description}
\item[$i \in I$] Buses.
\item[$j \in J$] Grid-connected devices.
\item[$k \in K$] Security contingencies.
\item[$t \in T$] Time intervals.
\item[$n \in N$] Reserve zones.
\item[$e \in E$] Miscellaneous constraints on operation of inedividual devices.
\end{description}

\subsubsection{Distinguished set elements}
\begin{description}
\item[$i_j \in I$] Connection bus of non-branch device $j$.
\item[$i^{fr}_j \in I$] From bus of branch $j$.
\item[$i^{to}_j \in I$] To bus of branch $j$.
\item[$j^{out}_k \in J$] The branch outaged in contingency $k$.
\end{description}

\subsubsection{Subsets}
\begin{description}
\item[$J^{pr} \subset J$] Producing devices, e.g. generators.
\item[$J^{cs} \subset J$] Consuming devices, e.g. loads.
\item[$J^{sh} \subset J$] Shunt devices.
\item[$J^{ac} \subset J$] AC branches.
\item[$J^{dc} \subset J$] DC branches.
\item[$J_i \subset J$] Non-branch devices at bus $i$, i.e. $J_i = \{ j \in J^{pr} \cup J^{cs} \cup J^{sh} : i_j = i \}$.
\item[$J^{fr}_i \subset J$] Branches from bus $i$, i.e. $J^{fr}_i = \{ j \in J^{ac} \cup J^{dc} : i^{fr}_j = i \}$.
\item[$J^{to}_i \subset J$] Branches to bus $i$, i.e. $J^{to}_i = \{ j \in J^{ac} \cup J^{dc} : i^{to}_j = i \}$.
\item[$J_k \subset J$] Devices remaining in service in contingency $k$, i.e. $J_k = J \setminus \{ j^{out}_k \}$.
\item[$J_n \subset J$] Producing or consuming devices contributing to reserve zone $n$.
\item[$T^{on}_j \subset T$] Must run intervals for producing or consuming device $j$.
\item[$T^{off}_j \subset T$] Forced outage intervals for producing or consuming device $j$.
\item[$T^{su}_{jt} \subset T$] Intervals prior to $t$ during which a shut down of producing or consuming device $j$ would preclude starting up in interval $t$ based on minimum downtime.
\item[$T^{sd}_{jt} \subset T$] Intervals prior to $t$ during which a start up of producing or consuming device $j$ would preclude shutting down in interval $t$ based on minimum uptime.
\item[$E^{su}_j \subset E$] Max startups constraints for producing or consuming device $j$.
\item[$E^{sd}_j \subset E$] Max shutdowns constraints for producing or consuming device $j$.
\item[$E^{en}_j \subset E$] Multi-interval total energy constraints on producing or consuming device $j$.
\end{description}

\subsubsection{Parameters}
\begin{description}
\item[$U^0_j$] Initial on-off status of producing or consuming device $j$, $1$ if on, $0$ if off, or status of AC branch $j$, $1$ if connected, $0$ if disconnected.
\item[$U^{max}_j$] Maximum number of steps that may be activated for shunt $j$.
\item[$U^{min}_j$] Minimum number of steps that may be activated for shunt $j$.
\item[$\tau^{max}_j$] Maximum value of winding ratio of AC branch $j$.
\item[$\tau^{min}_j$] Minimum value of winding ratio of AC branch $j$.
\item[$\theta^{max}_j$] Maximum value of phase difference of AC branch $j$.
\item[$\theta^{min}_j$] Minimum value of phase difference of AC branch $j$.
\item[$P^{max}_j$] Maximum real power flow on DC branch $j$ in either direction.
\item[$Q^{max,fr}_j$] Maximum reactive power flow into DC branch $j$ at from bus.
\item[$Q^{min,fr}_j$] Minimum reactive power flow into DC branch $j$ at from bus.
\item[$Q^{max,to}_j$] Maximum reactive power flow into DC branch $j$ at to bus.
\item[$Q^{min,to}_j$] Minimum reactive power flow into DC branch $j$ at to bus.
\item[$V^{max}_i$] Maximum voltage magnitude at bus $i$.
\item[$V^{max}_i$] Minimum voltage magnitude at bus $i$.
\item[$U^{su,max}_{je}$] Maximum number of startups for constraint $e$ of producing or consuming device $j$.
\item[$U^{sd,max}_{je}$] Maximum number of shutdowns for constraint $e$ of producing or consuming device $j$.
\item[$P^{max}_{jt}$] Maximum real power production (or consumption) of producing (or consuming) device $j$ if online in interval $t$.
\item[$P^{min}_{jt}$] Minimum real power production (or consumption) of producing (or consuming) device $j$ if online in interval $t$.
\item[$Q^{max}_{jt}$] Maximum reactive power production (or consumption) of producing (or consuming) device $j$ if online in interval $t$.
\item[$Q^{min}_{jt}$] Minimum reactive power production (or consumption) of producing (or consuming) device $j$ if online in interval $t$.
\item[$P^0_j$] Real power production (or consumption) of producing (or consuming) device $j$ prior to start of model horizon.
\item[$P^{ru}_j$] Ramp up rate of device $j$ while online.
\item[$P^{rd}_j$] Ramp down rate of device $j$ while online.
\item[$P^{su}_j$] Ramp up rate of device $j$ when starting up.
\item[$P^{sd}_j$] Ramp down rate of device $j$ when shutting down.
\item[$P^{rsv,max}_{jt}$] Maximum reserve provision of device $j$ in interval $t$.
\item[$Y^{sh}_j$] Step complex admittance to ground of shunt $j$.
\item[$Y^{sr}_j$] Complex admittance of the series element of AC branch $j$.
\item[$Y^{fr}_j$] Complex admittance to ground at the from bus of AC branch $j$.
\item[$Y^{to}_j$] Complex admittance to ground at the to bus of AC branch $j$.
\item[$B^{sr}_j$] Susceptance, i.e. imaginary part of complex admittance, of the series element of AC branch $j$.
\item[$\alpha_i$] Slack distribution factor of bus $i$ in the post-contingency lossless DC model, defined here as $\alpha_i = 1/|I|$.
\item[$D_t$] Duration of interval $t$, strictly positive.
\item[$C^{rsv}_j$] Cost of providing reserve of producing or consuming device$j$.
\item[$C^{su}_j$] Startup cost of producing or consuming device $j$.
\item[$C^{sd}_j$] Shutdown cost of producing or consuming device $j$.
\item[$C^{on}_j$] Fixed commitment (i.e. no load) cost of producing or consuming device $j$.
\item[$C^p$] Real power imbalance penalty coefficient.
\item[$C^q$] Reactive power imbalance penalty coefficient.
\item[$C^s$] Branch overload penalty coefficient.
\item[$C^{sw}$] Branch switching (i.e. reconnection or disconnection) cost.
\item[$C^{en}$] Penalty coefficient on violation of multi-interval total energy constraints.
\item[$C^{rsv}_n$] Reserve shortage penalty coefficient of reserve zone $n$.
\item[$S^{max}_j$] Apparent power flow limit on AC branch $j$.
\item[$S^{max,ctg}_j$] Post-contingency apparent power flow limit on AC branch $j$.
\item[$A^{en0}_{je}$] Constant term in multi-interval total energy constraint $e$ for producing or consuming device $j$.
\item[$A^{en}_{jet}$] Linear coefficient on power in interval $t$ in multi-interval total energy constraint $e$ for producing or consuming device $j$.
\end{description}

\subsubsection{Functions}
\begin{description}
\item[$S(w,w',Y,Y')$] Algebraic function defining complex power flow on an AC branch, specifically equal to the complex power flow
into the branch at the from bus if the branch has complex winding ratio equal to $1$, complex admittance $Y$ on the series element, complex admittance $Y'$ to ground at the from bus, when the complex voltages are $w$ and $w'$ at the from and to buses, and the branch is in service and closed.
\item[$C^{en}_{jt}(p)$] Concave (or convex) piecewise linear energy value (or cost) function of producing (or consuming) device $j$ in interval $t$, specifically equal to the value accruing to a consuming device (or cost incurred by a producing device)
assuming the power consumption (or production) is equal to $p$.
\end{description}

\subsubsection{Decision variables}
\begin{description}
\item[$v_{it}$] Voltage magnitude at bus $i$ in interval $t$.
\item[$\theta_{it}$] Voltage angle at bus $i$ in interval $t$.
\item[$w_{it}$] Complex voltage at bus $i$ in interval $t$.
\item[$p_{it}$] Real power imbalance (i.e. net excess withdrawal) at bus $i$ in interval $t$.
\item[$q_{it}$] Reactive power imbalance at bus $i$ in interval $t$.
\item[$s_{it}$] Complex power imbalance at bus $i$ in interval $t$.
\item[$p_{jt}$] Real power injection of a producing device $j$ or withdrawal of a consuming or shunt device $j$ in interval $t$.
\item[$q_{jt}$] Reactive power injection of a producing device $j$ or withdrawal of a consuming or shunt device $j$ in interval $t$.
\item[$s_{jt}$] Complex power injection of a producing device $j$ or withdrawal of a consuming or shunt device $j$ in interval $t$.
\item[$p^{fr}_{jt}$] Real power flow into AC or DC branch $j$ in interval $t$ at the from bus.
\item[$q^{fr}_{jt}$] Reactive power flow into AC or DC branch $j$ in interval $t$ at the from bus.
\item[$s^{fr}_{jt}$] Complex power flow into AC or DC branch $j$ in interval $t$ at the from bus.
\item[$\tau_{jt}$] Winding ratio of AC branch $j$ in interval $t$.
\item[$\phi_{jt}$] Phase difference of AC branch $j$ in interval $t$.
\item[$\nu_{jt}$] Complex winding ratio of AC branch $j$ in interval $t$.
\item[$u_{jt}$] Integer variable representing discrete settings of device $j$ in interval $t$. Represents commitment status for producing or consuming devices, number of steps activated for shunts, and connected (online or closed) or disconnected (offline or open) status for AC branches.
\item[$u^{su}_{jt}$] Binary variable indicating startup of producing or consuming device $j$ or reconnection of AC branch $j$ in interval $t$.
\item[$u^{sd}_{jt}$] Binary variable indicating shutdown of producing or consuming device $j$ or disconnection of AC branch $j$ in interval $t$.
\item[$p^{rsv}_{jt}$] Quantity of reserve provided by device $j$ in interval $t$. In this paper, only the synchronized (i.e. spinning) reserve product is formulated in detail, but the full formulation considers a suite of reserve products.
\item[$p^{rsv,req}_{nt}$] Reserve requirement in zone $n$ in interval $t$.
\item[$p^{rsv}_{nt}$] Reserve shortage in zone $n$ in interval $t$.
\item[$\theta_{itk}$] Post-contingency voltage angle at bus $i$ in interval $t$ under contingency $k$.
\item[$p_{jtk}$] Post-contingency real power flow on AC branch $j$ in interval $t$ under contingency $k$.
\item[$p^{sl}_t$] System-wide real power mismatch in the post-contingency lossless DC model in interval $t$.
\item[$z^{ms}$] Total market surplus objective, for maximization.
\item[$z^{base}$] Base case objective.
\item[$z^{ctg,worst}$] Total worst case post-contingency objective penalty.
\item[$z^{ctg,avg}$] Total average case post-contingency objective penalty.
\item[$z^{ctg,worst}_t$] Worst case post-contingency objective penalty for interval $t$.
\item[$z^{ctg,avg}_t$] Average case post-contingency objective penalty for interval $t$.
\item[$z^{time}_t$] Time-indexed component of the base case objective for interval $t$.
\item[$z^{en}_{je}$] Penalty on violation of multi-interval total energy constraint $e$ for producing or consuming device $j$.
\item[$z^{en}_{jt}$] Energy value accruing to consuming device $j$ or cost incurred by producing device $j$ in interval $t$.
\item[$z^{rsv}_{jt}$] Cost incurred by a producing or consuming device $j$ to provide reserve in interval $t$.
\item[$z^{su}_{jt}$] Startup cost incurred by producing or consuming device $j$ in interval $t$.
\item[$z^{su}_{jt}$] Shutdown cost incurred by producing or consuming device $j$ in interval $t$.
\item[$z^{on}_{jt}$] Fixed commitment (i.e. no load) cost incurred by producing or consuming device $j$ in interval $t$.
\item[$z^{sw}_{jt}$] Switching (i.e. reconnection or disconnection) cost incurred by AC branch $j$ in interval $t$.
\item[$z^s_{jt}$] Overload penalty incurred by AC branch $j$ in interval $t$.
\item[$z^{rsv}_{nt}$] Reserve shortage penalty incurred in reserve zone $n$ in interval $t$.
\item[$z^p_{it}$] Real power imbalance penalty incurred at bus $i$ in interval $t$.
\item[$z^q_{it}$] Reactive power imbalance penalty incurred at bus $i$ in interval $t$.
\item[$z^{ctg}_{tk}$] Post-contingency penalty in interval $t$ under contingency $k$.
\item[$z^s_{jtk}$] Post-contingency overload penalty incurred by AC branch $j$ in interval $t$ under contingency $k$.
\end{description}
\subsection{Sets and indices}
%
The network consists of buses $i \in I$ and grid-connected devices $j \in J$. The set of devices is partitioned into producing devices $J^{pr}$, consuming devices $J^{cs}$, shunt devices $J^{sh}$, AC branches $J^{ac}$, and DC branches $J^{dc}$. The time horizon consists of time intervals $t \in T$. We consider security contingencies $k \in K$. Reserves are modeled over reserve zones $n \in N$. A number of miscellaneous constraints on the operation of individual producing or consuming devices are indexed by $e \in E$.
\subsection{State variables and complex coordinates}
%
For each bus $i$ and time $t$, the voltage magnitude and angle and complex voltage are denoted by $v_{it}$, $\theta_{it}$, and $w_{it}$, and the net real, reactive, and complex power imbalance are denoted by $p_{it}$, $q_{it}$, and $s_{it}$. These variables are related by
\begin{gather}
    w_{it} = v_{it} \left( \cos{\theta_{it}} + \sqrt{-1} \sin{\theta_{it}} \right) \ \forall i \in I, t \in T \\
    s_{it} = p_{it} + \sqrt{-1} q_{it} \ \forall i \in I, t \in T
\end{gather}
The real, reactive, and complex power injection of a producer device or withdrawal of a consumer device or shunt $j$ at time $t$ at the connection bus $i_j$ are denoted by $p_{jt}$, $q_{jt}$, and $s_{jt}$ and are related by
\begin{gather}
    s_{jt} = p_{jt} + \sqrt{-1} q_{jt} \ \forall j \in J^{pr} \cup J^{cs} \cup J^{sh}, t \in T
\end{gather}
The real, reactive, and complex power flows into AC or DC branch $j$ at time $t$ at the from bus $i^{fr}_j$ (or the to bus $i^{to}_j$) are denoted by $p^{fr}$, $q^{fr}$, and $s^{fr}$ (or $p^{to}$, $q^{to}$, and $s^{to}$) and are related by:
\begin{gather}
    s^{fr}_{jt} = p^{fr}_{jt} + \sqrt{-1} q^{fr}_{jt} \ \forall j \in J^{ac}, t \in T \\
    s^{to}_{jt} = p^{to}_{jt} + \sqrt{-1} q^{to}_{jt} \ \forall j \in J^{ac}, t \in T
\end{gather}
The winding ratio, phase difference, and complex winding ratio of AC branch $j$ in time $t$ are denoted by $\tau_{jt}$, $\phi_{jt}$, and $\nu_{jt}$ and are related by:
\begin{equation}
    \nu_{jt} = \tau_{jt} \left( \cos{\phi_{jt}} + \sqrt{-1} \sin{\phi_{jt}} \right) \ \forall j \in J^{ac}, t \in T
\end{equation}
\subsection{Device settings and feasible regions}
Producing and consuming devices and shunts and AC branches $j$ have integer variables $u_{jt}$, representing discrete settings:
\begin{equation}
    u_{jt}, \in \{ \dots, -1, 0, 1, \dots \} \ \forall j \in J^{pr} \cup J^{cs} \cup J^{sh} \cup J^{ac}, t \in T
\end{equation}
For producing and consuming devices and AC branches, the $u_{jt}$ variables are binary, with the value $1$ representing the online or connected state and $0$ representing the offline or disconnected state:
\begin{equation}
    0 \leq u_{jt} \leq 1 \ \forall j \in J^{pr} \cup J^{cs} \cup J^{ac}, t \in T
\end{equation}
For producing and consuming devices, there may be must run and forced outage conditions:
\begin{gather}
    u_{jt} = 1 \ \forall j \in J^{pr} \cup J^{cs}, t \in T^{on}_j \\
    u_{jt} = 0 \ \forall j \in J^{pr} \cup J^{cs}, t \in T^{off}_j
\end{gather}
For producing and consuming devices and AC branches, we introduce additional binary variables $u^{su}_{jt}$ and $u^{sd}_{jt}$, with the value $1$ indicating startup (or shutdown) or in the case of AC branches connection (or disconnection):
\begin{gather}
    u^{su}_{jt}, u^{sd}_{jt} \in \{ 0, 1 \} \ \forall j \in J^{pr,cs,ac}, t \in T \\
    u_{jt} - 1_{t > 0} u_{j,t-1} - 1_{t = 0} U^0_j = u^{su}_{jt} - u^{sd}_{jt} \ \forall j \in J^{pr} \cup J^{cs} \cup J^{ac}, t \in T
\end{gather}
For shunts, the $u_{jt}$ variables take values between prescribed bounds and represent the number of steps that are activated
\begin{equation}
    U^{min}_j \leq u_{jt} \leq U^{max}_j \ \forall j \in J^{sh}, t \in T
\end{equation}
The winding ratio and phase difference of an AC branch is continuously adjustable within bounds, and these bounds are equal to $1$ and $0$ for transmission lines but may be different for tap changing or phase shifting transformers:
\begin{gather}
    \tau^{min}_j \leq \tau_{jt} \leq \tau^{max}_j \ \forall j \in J^{ac}, t \in T \\
    \phi^{min}_j \leq \phi_{jt} \leq \phi^{max}_j \ \forall j \in J^{ac}, t \in T
\end{gather}
DC branches have bounded flows and $0$ real power losses:
\begin{gather}
    -P^{max}_j \leq p^{fr}_{jt} \leq P^{max}_j \ \forall j \in J^{dc}, t \in T \\
    -P^{max}_j \leq p^{to}_{jt} \leq P^{max}_j \ \forall j \in J^{dc}, t \in T \\
    Q^{min,fr}_j \leq q^{fr}_{jt} \leq Q^{max}_j \ \forall j \in J^{dc}, t \in T \\
    Q^{min,to}_j \leq q^{to}_{jt} \leq Q^{max}_j \ \forall j \in J^{dc}, t \in T \\
    p^{fr}_{jt} + p^{to}_{jt} = 0 \ \forall j \in J^{dc}, t \in T
\end{gather}
Bus voltage magnitudes are bounded:
\begin{equation}
    V^{min}_i \leq v_{it} \leq V^{max}_i \ \forall i \in I, t \in T
\end{equation}
\subsection{Producing and consuming device scheduling constraints}
Producing and consuming devices have scheduling constraints restricting the number of startups or shutdowns across sets of intervals:
\begin{gather}
    \sum_{t \in T^{su}_e} u^{su}_{jt} \leq U^{su,max}_{je} \ \forall j \in J^{pr} \cup J^{cs}, e \in E^{su}_j \\
    \sum_{t \in T^{sd}_e} u^{sd}_{jt} \leq U^{sd,max}_{je} \ \forall j \in J^{pr} \cup J^{cs}, e \in E^{sd}_j
\end{gather}
Minimum uptime (and downtime) constraints preclude startup (or shutdown) if a device has shut down (or started up) within a prescribed set of prior intervals:
\begin{gather}
    u^{su}_{jt} \leq 1 - \sum_{t' \in T^{dn}_{jt}} u^{sd}_{jt'} \ \forall j \in J^{pr} \cup J^{cs}, t \in T \\
    u^{sd}_{jt} \leq 1 - \sum_{t' \in T^{up}_{jt}} u^{su}_{jt'} \ \forall j \in J^{pr} \cup J^{cs}, t \in T
\end{gather}
\subsection{Producing and consuming device real and reactive power and reserves}
The real and reactive power of producing and consuming devices $j$ is bounded by time-varying limits:
\begin{gather}
    P^{min}_{jt} u_{jt} \leq p_{jt} \leq P^{max}_{jt} u_{jt} \ \forall j \in J^{pr} \cup J^{cs}, t \in T \\
    Q^{min}_{jt} u_{jt} \leq q_{jt} \leq Q^{max}_{jt} u_{jt} \ \forall j \in J^{pr} \cup J^{cs}, t \in T
\end{gather}
Ramp rates constrain the change in real power over time:
\begin{multline}
    -D_t \left( P^{rd}_j u_{jt} + P^{sd}_j u^{su}_{jt} \right)
\leq p_{jt} - 1_{t > 0} p_{j,t-1} - 1_{t = 0} P^0_j \\ \leq D_t \left( P^{ru}_j ( u_{jt} - u^{su}_{jt} ) + P^{su}_j u^{su}_{jt} \right) \ \forall j \in J^{pr} \cup J^{cs}, t \in T
\end{multline}
%
The full formulation considers a suite of reserve products, including regulation up and down, synchronized (i.e. spinning) and non-synchronized reserve, ramping reserve up and down, and reactive power reserve up and down. Both producing and consuming devices are modeled as able to provide these products, within technical capabilities. For this paper, we formulate in detail just the spinning reserve product, with a variable $p^{rsv}_{jt}$ denoting the quantity of reserve provided by device $j$ in interval $t$:
\begin{gather}
    0 \leq p^{rsv}_{jt} \leq P^{rsv,max}_{jt} u_{jt} \ \forall j \in J^{pr} \cup J^{cs}, t \in T \\
    p_{jt} + p^{rsv}_{jt} \leq P^{max}_{jt} u_{jt} \ \forall j \in J^{pr}, t \in T \\
    p_{jt} - p^{rsv}_{jt} \geq P^{min}_{jt} u_{jt} \ \forall j \in J^{cs}, t \in T
\end{gather}
Other aspects of producing and consuming device operations modeled in the full formulation 
but not formulated here include: reactive power capability constraints limiting both real and reactive power injection to a trapezoidal set or a line segment for constant power factor, multi-interval startup and shutdown real power trajectories, and multi-interval total energy constraints to model storage and time-shiftable load. 
%
%
\subsection{AC power flow and balance}
Power flows into shunts $j$ are determined by the number of activated steps and the bus voltage as well as the shunt step admittance $Y^{sh}_j$
\begin{gather}
    s_{jt} = Y^{sh*}_j u_{jt} v_{it}^2 \ \forall j \in J^{sh}, t \in T, i = i_j
\end{gather}
%
%
%
%
%
%
The power flows on an AC branch $j$ are equal to $0$ if the branch is disconnected and are otherwise a function of the voltages at the from and to buses and the complex winding ratio on the branch as well as the branch admittances $Y^{fr}_j$, $Y^{to}_j$, and $Y^{sr}_j$ at the from and to buses and on the series element of the branch:
\begin{gather}
    s^{fr}_{jt} = u_{jt} S(w_{it} / \nu_{jt}, w_{i't}, Y^{sr}_j, Y^{fr}_j) \ \forall j \in J^{ac}, t \in T, i = i^{fr}_j, i' = i^{to}_j \\
    s^{to}_{jt} = u_{jt} S(w_{i't}, w_{it} / \nu_{jt}, Y^{sr}_j, Y^{to}_j) \ \forall j \in J^{ac}, t \in T, i = i^{fr}_j, i' = i^{to}_j
\end{gather}
The function defining flows on an AC branch is
\begin{equation}
    S(w,w',Y,Y') = Y'^* w w^* + Y^* w (w - w')^*
\end{equation}
Power balance:
\begin{equation}
    \sum_{j \in (J^{cs} \cup J^{sh}) \cap J_i} s_{jt}
    - \sum_{j \in J^{pr} \cap J_i} s_{jt}
    + \sum_{j \in J^{fr}_i} s^{fr}_{jt}
    + \sum_{j \in J^{to}_i} s^{to}_{jt}
    = s_{it}
    \ \forall i \in I, t \in T
\end{equation}
%
\subsection{Security constraints}
Security constraints require that the planned dispatch be secure to any of a set of contingencies $k \in K$, each defined by the unplanned outage of a single branch $j^{out}_k$, so that the set of remaining devices is $J_k = J \setminus \{ j^{out}_k \}$.
Post-contingency AC branch real power flows $p_{jtk}$ follow a DC flow model, with variables $\theta_{itk}$ representing the post-contingency bus voltage angles, and variables for AC branch on-off status and phase difference fixed to their pre-contingency values, using just the imaginary part $B^{sr}_j$ of the branch admittance $Y^{sr}_j$:
\begin{equation}
    p_{jtk} = -u_{jt} B^{sr}_j \left( \theta_{itk} - \theta_{i'tk} - \phi_{jt} \right) \ \forall t \in T, k \in T, j \in J_k, i = i^{fr}_j, i' = i^{to}_j
\end{equation}
Real power injections and withdrawals of non-branch devices and remaining DC lines are fixed to their pre-contingency values, and real power balance is enforced at each bus. Since the pre-contingency real power values reflect an AC model that includes losses, and the post-contingency AC branch flows follow a lossless DC model, there may be a nonzero system-wide real power mismatch that is represented by a variable $p^{sl}_t$:
\begin{equation}
    p^{sl}_t = \sum_{j \in J^{pr}} p_{jt} - \sum_{j \in J^{cs}} p_{jt} - \sum_{j \in J^{sh}} p_{jt} \ \forall t \in T
\end{equation}
The system slack is distributed across the buses $i$ in proportion to fixed slack distribution coefficients, which we define uniformly as $\alpha_i = 1 / |I|$. Then the post-contingency power balance constraints are
\begin{multline}
    \sum_{j \in J^{ac} \cap J^{fr}_i \cap J_k} p_{jtk}
    - \sum_{j \in J^{ac} \cap J^{to}_i \cap J_k} p_{jtk}
    + \alpha_i p^{sl}_t
    = \sum_{j \in J^{pr} \cap J_i} p_{jt} \\
    - \sum_{j \in (J^{cs} \cup J^{sh}) \cap J_i} p_{jt}
    - \sum_{j \in J^{dc} \cap J^{fr}_i \cap J_k} p^{fr}_{jt}
    - \sum_{j \in J^{dc} \cap J^{to}_i \cap J_k} p^{to}_{jt}
    \\ \forall i \in I, t \in T, k \in K
\end{multline}
\subsection{Connectivity constraints}
We require that the network consisting of all buses and closed branches be connected in every time period $t$, in the base case and in every contingency $k$. In other words, for every interval $t$, every pair of buses should be connected via a path of branches in the set $J^{dc} \cup \{ j \in J^{ac} : u_{jt} = 1 \}$. Furthermore, for every interval $t$ and every contingency $k$, every pair of buses should be connected via a path of branches in the set $(J^{dc} \cap J_k) \cup \{ j \in J^{ac} \cap J_k : u_{jt} = 1 \}$. This requirement is needed for realism, but it has the further effect of ensuring that, given the base case solution, there is a unique post-contingency solution, up to a network-wide shift of voltage angles. We do not formulate this requirement as an algebraic constraint. We do note that it is expressed in terms of the values of the finitely many binary variables $u_{jt}$ for $j \in J^{ac}$ and $t \in T$ and therefore describes a finite set of points in the space of these variables.
\subsection{Market Surplus Objective}
The total market surplus objective $z^{ms}$ for maximization is the sum of a base case objective $z^{base}$ minus a worst case contingency penalty $z^{ctg,worst}$ and an average case contingency penalty $z^{ctg,avg}$:
\begin{equation}
    \label{eq:z_ms}
    z^{ms} = z^{base} - z^{ctg,worst} - z^{ctg,avg}
\end{equation}
The worst case and average case contingency penalties are the sum of time-indexed worst case and average case contingency penalties $z^{ctg,worst}_t$ and $z^{ctg,avg}_t$ over times $t$:
\begin{gather}
    \label{eq:z_ctg_worst}
    z^{ctg,worst} = \sum_{t \in T} z^{ctg,worst}_t \\
    \label{eq:z_ctg_avg}
    z^{ctg,avg} = \sum_{t \in T} z^{ctg,avg}_t
\end{gather}
The base case objective is the sum of the time-indexed base case objectives $z^{time}_t$ for times $t$ minus penalties $z^{en}_{je}$ for multi-interval energy constraints $e \in E^{en}_j$ on producing and consuming devices $j \in J^{pr} \cup J^{cs}$ that cannot be indexed to time:
\begin{equation}
    \label{eq:z_base}
    z^{base} = \sum_{t \in T} z^{time}_t - \sum_{j \in J^{pr} \cup J^{cs}, e \in E^{en}_j} z^{en}_{je}
\end{equation}
The time-indexed base case objective is the sum of appropriately signed terms representing the values accrued and costs and penalties incurred for consumption and production of energy, provision of reserves, unit scheduling, topology switching, transmission overloading, reserve shortage, and real and reactive power imbalance:
\begin{multline}
    z^{time}_t = \sum_{j \in J^{cs}} z^{en}_{jt} - \sum_{j \in J^{pr}} z^{en}_{jt} - \sum_{j \in J^{pr} \cup J^{cs}} \left( z^{rsv}_{jt} + z^{su}_{jt} + z^{su}_{jt} + z^{on}_{jt} \right) \\ - \sum_{j \in J^{ac}} \left( z^{sw}_{jt} + z^s_{jt} \right) - \sum_{n \in N} z^{rsv}_{nt} - \sum_{i \in I} \left( z^p_{it} + z^q_{it} \right)
\end{multline}
Specifically, in time $t$, $z^{en}_{jt}$ is the value accruing to consuming device $j \in J^{cs}$ from consumption of energy or the cost incurred by producing device $j \in J^{pr}$ from production of energy; $z^{rsv}_{jt}$, $z^{su}_{jt}$, $z^{sd}_{jt}$, and $z^{on}_{jt}$ are the cost of providing reserves and startup, shutdown, and no-load costs incurred by producing or consuming device $j \in J^{pr} \cup J^{cs}$; $z^{sw}_{jt}$ and $z^s_{jt}$ are the cost of topology switching and the penalty of transmission limit overload incurred by AC branch $j \in J^{ac}$; $z^{rsv}_{nt}$ is the penalty incurred by reserve shortage in reserve zone $n \in N$; and $z^p_{it}$ and $z^q_{it}$ are the penalties incurred by real and reactive power imbalance at bus $i \in I$.
For each interval, the time-indexed worst case and average case contingency penalties are defined as the maximum and the average over contingencies $k \in K$ of contingency-indexed penalties $z^{ctg}_{tk}$:
\begin{gather}
    z^{ctg,worst}_t = \max_{k \in K} z^{ctg}_{tk} \ \forall t \in T \\
    z^{ctg,avg}_t = 1 / |K| \sum_{k \in K} z^{ctg}_{tk} \ \forall t \in T
\end{gather}
The contingency-indexed penalties are the sum over AC branches $j$ of post-contingency transmission limit overload penalties $z^s_{jtk}$:
\begin{equation}
    z^{ctg}_{tk} = \sum_{j \in J^{ac} \cap J_k} z^s_{jtk} \ \forall t \in T, k \in K
\end{equation}
\subsection{Components of objective terms}
The value accruing to a consuming device (or the cost incurred by a producing device) $j$ from the consumption (or production) $p_{jt}$ of power during time $t$ is modeled with a concave (or convex) piecewise linear energy value (or cost) function $C^{en}_{jt}$:
\begin{equation}
    z^{en}_{jt} = D_t C^{en}_{jt} (p_{jt}) \ \forall j \in J^{pr} \cup J^{cs}, t \in T
\end{equation}
Here $D_t$ is the duration of time interval $t$, which may differ across intervals. Throughout the formulation, a factor of $D_t$ is applied to time-intensive (i.e. per hour) quantities to transform them into absolute quantities.
The cost to a producing or consuming device $j$ of providing reserves 
is a linear function of the quantity of reserve provided: 
\begin{equation}
    z^{rsv}_{jt} = D_t C^{rsv}_j p^{rsv}_{jt} \ \forall j \in J^{pr} \cup J^{cs}, t \in T
\end{equation}
%
Startup, shutdown, and no-load cost are formulated in a simple fashion here, but the full formulation 
includes downtime-dependent startup costs as well:
\begin{gather}
    z^{su}_{jt} = C^{su}_j \max \left(
    0, u_{jt} - 1_{t > 0} u_{j,t-1} - 1_{t = 0} U^0_j
    \right) \ \forall j \in J^{pr} \cup J^{cs}, t \in T \\
    z^{sd}_{jt} = C^{sd}_j \max \left(
    0, 1_{t > 0} u_{j,t-1} + 1_{t = 0} U^0_j - u_{jt}
    \right) \ \forall j \in J^{pr} \cup J^{cs}, t \in T \\
    z^{on}_{jt} = D_t C^{on}_j u_{jt} \ \forall j \in J^{pr} \cup J^{cs}, t \in T
\end{gather}
%
%
%
Bus real and reactive power imbalance penalties are multiples of the absolute values of the corresponding imbalances:
\begin{gather}
    z^p_{it} = D_t C^p |p_{it}| \ \forall i \in I, t \in T \\
    z^q_{it} = D_t C^q |q_{it}| \ \forall i \in I, t \in T
\end{gather}
Penalties are applied to zonal reserve shortages $p^{rsv}_{nt}$,
and for the generic reserve product modeled here, the reserve requirement $p^{rsv,req}_{nt}$ for reserve zone $n$ is a multiple $\sigma_n$  of the largest energy dispatch value over producing devices in zone $n$, where $J_n$ is the set of devices in zone $n$:
\begin{gather}
    z^{rsv}_{nt} = D_t C^{rsv}_n p^{rsv}_{nt} \ \forall n \in N, t \in T \\
    p^{rsv}_{nt} = \max \left( 0, p^{rsv,req}_{nt} - \sum_{j \in (J^{pr} \cup J^{cs}) \cap J_n} p^{rsv}_{jt} \right) \ \forall n \in N, t \in T \\
    p^{rsv,req}_{nt} = \sigma_n \max_{j \in J^{pr} \cap J_n} p_{jt} \ \forall n \in N, t \in T
\end{gather}
In the full formulation, the requirements of the various reserve products include endogenous requirements (i.e. depending on the dispatch, as with the generic product modeled here) and exogeneous requirements (i.e. a fixed value for each zone and time period). Some of the endogenous requirements (e.g. spinning reserve) depend on the maximum device energy dispatch, while others (e.g. regulation) depend on the total over all devices.
Multi-interval energy limits are formulated as soft constraints with penalties on violations:
\begin{equation}
    z^{en}_{je} = C^{en} \max \left( 0, A^{en0}_{je} + \sum_{t \in T} A^{en}_{jet} p_{jt} \right) \ \forall j \in J^{pr} \cup J^{cs}, e \in E^{en}_j
\end{equation}
Topology switching cost on AC branches
\begin{equation}
    z^{sw}_{jt} = D_t C^{sw} \left|
    u_{jt} - 1_{t > 0} u_{j,t-1} - 1_{t = 0} U^0_j
    \right| \ \forall j \in J^{ac}, j \in T
\end{equation}
AC branch flow limits enforced as soft constraints with penalties on violations
\begin{equation}
    z^s_{jt} = D_t C^s
    \max \left( 0, \max \left( |s^{fr}_{jt}|, |s^{to}_{jt}| \right) - S^{max}_j \right)
    \ \forall j \in J^{ac}, t \in T
\end{equation}
Post-contingency AC branch flow limits are similarly enforced as soft constraints, with reactive power flows assumed to be as in the base case and real power flows $p_{jtk}$ defined by the post-contingency DC flow model:
\begin{multline}
    z^s_{jtk} = D_t C^s \max \left( 0, \left( p_{jtk}^2 + \max ( q^{fr}_{jt}, q^{to}_{jt} )^2 \right)^{1/2} - S^{max,ctg}_j \right) \\ \forall k \in K, j \in J^{ac} \cap J_k, t \in T
\end{multline}
\section{Traditional Approaches}
\label{sec:sol_methods_traditional}

{Unit commitment problem, UCP, is a day-ahead market clearing system operation practice that determines the discrete on/off decisions and continuous dispatches of participating generation units. The problem is known as mixed integer nonlinear programming, or MINLP. For decades, linearized DC power flow models have been used to reduce computational complexity when describing how nodal voltages affect power flows across networks. The problem for a large-scale network is frequently extremely complex. There have been two approaches that have received the most attention: dynamic programming and mixed integer linear programming, or MILP.

The sequential nature of decision-making inherent in the UCP lends itself particularly well to Dynamic Programming (DP). It divides the problem into stages and solves each one optimally. If all stages are solved correctly, this methodology ensures that the optimal schedule is identified. Furthermore, DP can handle a wide range of constraints and nonlinearities within a stage, which is common in power generation problems. Dynamic programming identifies the UC schedule to intermingle the temporal stages if the OPFs considered in each stage are all independent. The major disadvantage of DP is the "curse of dimensionality." The computational burden grows exponentially as the problem size grows, making it impractical for large-scale power systems. This limitation is significant in real-world applications with a large number of units and time periods. Furthermore, DP implementation can be complex and may necessitate significant computational resources. However, because the complexity is considered at each stage, this DP approach can employ the highly complex nature of power flow.
 
Due to advancements in MILP solvers and computational power, Mixed Integer Linear Programming (MILP) has become increasingly popular for solving the UCP. MILP is effective at dealing with large-scale problems and strikes an appropriate balance between solution quality and computational efficiency. The UCP is formulated as a linear programming problem with integer constraints, making it tractable for large datasets. MILP solvers use sophisticated algorithms to efficiently explore the solution space, which frequently results in the discovery of good solutions in reasonable time frames. The disadvantage of MILP is its linear approximation. MILP may require additional constraints and variables to approximate nonlinearities in real-world power systems, which can complicate the model and potentially impact the accuracy of the solution. Furthermore, while MILP is more efficient than DP at handling large problems, there is still a trade-off between solution accuracy and computational time, particularly for very large and complex power systems.
Although DP is very accurate and best suited for smaller or more divided problems, its computational limitations prevent it from being applied to large-scale scenarios. While linear approximations cause some precision loss, MILP provides a more realistic approach for large-scale problems by balancing computational viability and solution accuracy. The particular needs and size of the UCP at hand are often determining factors when choosing between DP and MILP.
 
The problem proposed in GO Ch. 3 differs from the traditional UCP in that it incorporates the nonconvex full  alternating current (AC) power flow model. The addition of the full AC power flow model to the UCP significantly increases the problem's complexity. For realistic modeling of power system operations, the AC power flow equations, which include nonlinear relationships, must be accurately represented. The nonlinear AC power flow model is written in Mixed Integer Nonlinear Programming (MINLP). Nonlinear equations can be directly included in the model using MINLP. However, solving MINLP problems is more computationally difficult than solving MILP problems, especially for large-scale power systems. The combination of integer decision variables and nonlinear equations adds to the complexity. Several approaches can be considered when transitioning from a linear power flow model to a full AC model within the MILP framework: relaxation techniques, decomposition techniques, extended DP, and hybrid approaches.
 
Linearization of AC power flow is the process of linearizing the equations of AC power flow around a specific operating point. This method converts nonlinear equations into linear equations that can be easily integrated into the MILP framework. However, this linearization can result in inaccuracies, particularly when large deviations from the chosen operating point are present. Relaxation Relaxing the nonlinear AC power flow equations can be accomplished using techniques such as second-order cone programming (SOCP) or semidefinite programming (SDP). These methods use convex constraints to approximate nonlinear constraints, making the problem more tractable while maintaining a higher degree of accuracy than simple linearization. This method strikes a reasonable balance between computational feasibility and accuracy.
 
Decomposition Techniques entails employing decomposition techniques to solve the UCP and AC power flow problems iteratively. The UCP problem can be solved using a simplified power flow model, and the generated generation schedule can then be evaluated using the full AC power flow model. If the AC power flow constraints are violated, the UCP solution is adjusted, and the process is repeated until a viable solution is found. According to the problem formulation described in the previous section, the MINLP incorporates soft constraints rather than hard constraints using linear penalty functions, and the UCP solution always returns a feasible solution. The formulation of the problem may encourage this approach, and as a result, we observed that many competitors use these techniques.
 
Combining methods can also be beneficial. Critical time periods or scenarios, for example, may require a detailed AC model, whereas less critical periods may require a simplified model. This method's primary advantage is its balance of accuracy and computational efficiency. It ensures that the UCP is solved precisely enough during critical periods while maintaining overall computational tractability. However, the added complexity of managing two different modeling paradigms and ensuring consistency and accuracy in the transition between them is a disadvantage. There's also a chance that the simplified model will overlook certain system nuances during non-critical periods, resulting in suboptimal or less robust decision-making. Therefore, the effectiveness of this approach is heavily reliant on correctly identifying which periods necessitate detailed modeling and which do not, which necessitates a thorough understanding of the system's dynamics and operational characteristics.
 
Each approach has trade-offs between computational demand, solution accuracy, and the ability to deal with the complexities of real-world power systems. The method chosen is determined by the specific requirements of the power system under consideration, the available computational resources, and the level of accuracy required.}

\section{Solution Approaches of the Competition Entrants}
\label{sec:sol_methods_entrants}

The major accomplishments of the GO competitors in solving the complex integration of Unit Commitment (UC) and AC Optimal Power Flow (AC OPF) have predominantly followed a decomposition approach. Various teams, while employing diverse methodologies, have commonly decomposed the problem into an initial UC phase coupled with a simplified OPF to ascertain the UC decisions. A critical observation across these methodologies is the importance of how the continuous variables, determined by the UC decisions, are managed in the subsequent OPF phase. This integrated analysis delves into the comparative efficiency of these approaches, focusing on the utilization of continuous variables post-UC decisions.

In sections \ref{sec:approach_decomp} - \ref{sec:approach_variations}, we present general observations on the approaches of the competition entrants. Then in Section \ref{sec:approach_specific}, we provide individual descriptions of these approaches based on information submitted to the competition organizers by the entrants themselves.

\subsection{Common Decomposition Strategy: UC Followed by Simplified OPF}

\label{sec:approach_decomp}
A recurring theme in the methodologies is the bifurcation of the problem into two distinct phases. Initially, the UC problem is solved to establish the operational decisions for generation units. This phase generally involves simplifying assumptions and heuristic techniques to ensure computational tractability. The teams adopted various strategies, such as iterative rounding, linear approximation, and relaxation techniques, to efficiently navigate through the complex landscape of integer and continuous variables in the UC problem.

\subsection{The Role of Continuous Variables in Post-UC OPF}
\label{sec:approach_cts}
{Post-UC decisions, the focus shifts to the OPF phase, where the problem becomes nearly decomposable across time steps, allowing for the application of heuristic approaches. A key differentiator in performance among the teams was noted in the handling of continuous variables determined from the UC phase. Teams that strategically accepted and integrated these continuous variables into the OPF phase demonstrated better performance. This approach effectively bridges the initial UC decisions with the more detailed OPF analysis, ensuring that the continuous variables such as power flow and voltage levels are realistically aligned with the discrete decisions made during the UC phase.}

\subsection{Methodological Variations and Performance}
\label{sec:approach_variations}
Each team, while conforming to the general structure of problem decomposition, exhibited variations in the integration and handling of continuous variables:

Performer A's methodology involved solving the ACOPF twice, each time with different constraints and considerations regarding reserves and voltages. This dual-phase approach allowed for a refined handling of continuous variables, albeit with the trade-off of increased computational complexity.

Performer B addressed the challenge of integrating AC losses into the UC model. By adjusting the continuous power flow variables in the UC phase based on previous operational points, they effectively bridged the gap between the discrete and continuous aspects of the problem.

Performer C and others employed relaxation and approximation techniques in the UC phase, focusing on achieving a balance between computational efficiency and the realistic representation of continuous variables in the OPF phase.

\subsection{Concluding Remarks}

{The comparative analysis of these methodologies underscores a critical aspect of power system optimization: the effective management of continuous variables post-UC decisions. This approach not only ensures coherence between the UC and OPF phases but also significantly enhances the overall performance and realism of the solution. The diversity in methods reflects the complexity and multi-dimensional nature of the UC and AC OPF integration challenge, highlighting the need for innovative, adaptable, and efficient solutions in this evolving field of power system optimization.}

\subsection{Individual Solution Approaches of the Competition Entrants}
\label{sec:approach_specific}

In sections \ref{sec:a} - \ref{sec:f}, we provide descriptions of the individual solution approaches of the competition entrants, based on information submitted to the competition organizers by the entrants themselves. These descriptions contain as much detail as the entrants are willing to make public at this time, and we are obligated to maintain this degree of confidentiality. In particular, the teams are identified anonymously and thus cannot be mapped to their public identifying information, which we do supply here. For convenient reference, we list here the key characteristics of the different approaches, including the main software, ideas, and techniques used:
\begin{itemize}
\item Performer A
\begin{itemize}
\item Julia, JuMP, Gurobi, Ipopt
\item UC over all periods, highly simplified
\item Simplification, relaxation, restriction, and approximation
\item ACOPF for each period, ignore reserve and solve, then fix voltage and flows and solve again for reserve.
\end{itemize}
\item Performer B
\begin{itemize}
\item Continuous relaxation of integer variables yielding a multiperiod ACOPF
\item Novel iterative batch rounding heuristic: solve multiperiod ACOPF, round and fix discrete variables in carefully chosen batches
\item Novel linearization of AC power flow based on losses at previous iteration
\end{itemize}
\item Performer C
\begin{itemize}
\item Python, Gurobi through Pyomo for MIP solver for UC model, Ipopt through Cyipopt for NLP solver for ACOPF model
\item McCormick envelope for tight relaxation of bilinear terms in quadratic constraints added to UC model
\item Novel heuristic seeking variable values that can be maintained across all time periods, reduces multiperiod ACOPF to single period
\end{itemize}
\item Performer D
\begin{itemize}
\item C, Lagrangian decomposition for UC, custom interior point NLP solver for ACOPF model, Paradiso for linear algebra
\item Decouple ACOPF across time periods by setting bounds to enforce ramping
\item AC violations fed back to UC through linear sensitivity factors
\item HPC computation of shift factors for power flow solution
\end{itemize}
\item Performer E
\begin{itemize}
\item Python, Ampl, Knitro, no discrete optimization solver
\item Handle discrete variables by fixing some, relaxing the others, solve NLP to low accuracy
\item Intelligent round and fix heuristic: conservative rounding prioritizing feasibility
\item Linear approximation of power balance equations
\end{itemize}
\item Performer F
\begin{itemize}
\item Gurobi for MIP solver for UC model, Ipopt with MA97 for NLP solver for ACOPF model
\item UC using a linear approximation of active power flow on the full time horizon
\item ACOPF on a single time step
\item UC and ACOPF linked via ADMM algorithm
\item Voltage violations handled by lazy constraints in the UC solver
\item Rank 1 update ofPTDF for fast contingency analysis rather than time consuming exact evaluation of objective at each iteration
\item Write out solution at each iteration if objective appears to be improved
\end{itemize}
\end{itemize}

Because of our commitment to the confidentiality of the information supplied to us by the competition entrants, we are unable to directly reference the solution approaches. Some of the entrants have published articles or given public presentations on their approaches, and these public materials are listed on the GO Competition web pages on citations \cite{go_competition_citations_web_page_misc} and on news \cite{go_competition_news_web_page_misc}. We encourage the reader to contact individual solver team members for further details on their approaches, and to facilitate this, we provide the following list of the solver teams, with their lead organizations and individual members:
\begin{itemize}
    \item ARPA-e Benchmark. Los Alamos National Laboratory. Carleton Coffrin, Robert Parker.
\item Argonauts. Argonne National Laboratory. Kibaek Kim, Michel Schanen, Anirudh Subramanyam, Weiqi Zhang, Francois Pacaud.
\item Artelys\_Columbia. Artelys. Richard A Waltz, Daniel Bienstock.
\item Electric-Stampede. University of Texas - Austin. Javad Mohammadi, Kyri Baker, Constance Crozier, Hussein Sharadga.
\item Enline. Enline (Brazil/Portugal). Leonardo Bonatto, Felipe da Silveira, Matheus Barbosa de Souza.
\item Gatorgar. University of Florida. Yongpei Guan, Lei Fan, Mingze Li, Weihang Ren.
\item GOT-BSI-OPF. Global Optimal Technology, Inc.. Hsiao-Dong Chiang, Pat Causgrove, Bin Wang, Lin Zeng.
\item GravityX. Individual. Hassan Hijazi.
\item LLGoMax. Lawrence Livermore National Laboratory. Ignacio Aravena Solis, Rod Frowd, Shmuel Oren, Alex Papalexopoulos.
\item Occam's razor. Individual. Wenyuan Tang.
\item PACE. IncSys. Robin Podmore, Roozbeh Khodadadeh, Chris Mosier.
\item PGWOpt. University of Pittsburgh. Masoud Barati, Nina Fatehi, Sattiago Grijalva, Masoud H. Nazari, Samuel Talkington, Jorge Fernandez.
\item Powersense. Individual. Sarah Sadat, Sayed Sadat.
\item quasiGrad. Individual. Samuel Chevalier.
\item The Blackouts. University of Tennessee - Knoxville. James Ostrowski, Marzieh Bakhshi, Christopher Ginart, William Hart, Bernard Knueven, Jonathan Schrock, Jean-Paul Watson.
\item TIM-GO. Massachusetts Institute of Technology. Xu Sun, Matthew Brun, Xin Chen, Dirk Lauinger, Thomas Lee.
\item xtellix. Xtellix, Inc. Mark Amo-Boateng.
\item YongOptimization. Mississippi State University. Yong Fu, Lin Gong, Yehong Peng, Fasiha Zainab.
\end{itemize}

\subsubsection{Performer A}
\label{sec:a}

This team did the entire project in about six weeks. Thus a series of substantial simplifications were made with the goal of achieving a limited quality of the approach. Quick decisions on the tools were made by choosing Julia with JuMP, Gurobi for discrete and IPOPT for continuous optimization. Intense experimentation lead to the ultimately submitted program.

To achieve a feasible solution in reasonable time the following major simplifications are made. UC is done for the entire time horizon and in each time period ACOPF is solved twice. Each subproblem is simplified and the feasible region modified through relaxation, restriction, and approximation. The key is to find adequate alterations and (if needed) projections onto the original feasible region.

The UC part of the algorithm was made simple or even very simple, One reason for solving the ACOPF part twice is to be able to ignore reserves in the first round and fix voltages and flows. A second reason is to assure numerical stability and precision. The resulting method is simple to understand and transparent, achieving the desired moderate solution quality.

This general algorithm, however, cannot work well for all datasets. So, several
parameters such as the choice of UC, whether reserves are used in the first ACOPF round, and several others were customized for different datasets.

In a future continuation the algorithm would be made more sophisticated through
a rolling horizon approach and a spatial decomposition. Parallel computation will
also be applied to accelerate the computation.

\subsubsection{Performer B}
\label{sec:b}

To claim mastery of a given problem, one can derive scientific formulas, write concise
mathematical theorems, and publish articles in prestigious venues, but when push comes to
shove, the reality of solving a real-world problem can hit hard.
The importance of having a well-defined and transparent evaluation procedure with realistic
datasets cannot be stressed enough. With different teams competing for the best approach, the
pressure of performing well and getting to the top of the leaderboard can be very motivating.
Most importantly, under this setup, there is zero tolerance for false claims or fake contributions,
if your approach is not correct or if it does not scale, it will soon be exposed, your leaderboard
score will reflect it.

One main challenge in the GO competition consisted in dealing with integer variables.
Performer B adopted a fast heuristic outlined below.
The Iterative Batch Rounding (IBR) heuristic starts by solving the continuous relaxation of the
given MINLP. It then loops over all variable batches, rounding and fixing the discrete variables
appearing in each batch separately before rerunning the reduced MINLP. The heuristic iterates
till all batches have been covered. The ordering of the batches, as well as the custom rounding
method play an important role.
Another major challenge consisted in capturing nonlinear constraints such as AC losses in the
Unit Commitment (UC) model:
Using lossless models often led to the under-commitment of generators, consequently leading
to large penalties due to soft constraint violations due to nodal imbalances.
The approach adopted by Performer B consists of using the AC losses obtained by fixing
the UC binaries to the previous operating point solution, solving the resulting ACOPF, computing
the active and reactive loss on every edge, then adding the linear constraints enforcing linear
loss constraints on active and reactive power flow variables appearing in the UC MIP
formulation. This approach can be iteratively refined by solving an ACOPF after each MIP and
recomputing the loss terms until a fixed point is reached or the time limit is exceeded.

Combining spatial decomposition with temporal decomposition and mixing coarse-grained with
fine-grained models is a future research direction. Combining machine learning with
mathematical optimization is also another promising research avenue.

\subsubsection{Performer C}
\label{sec:c}
  
{This team also used the split DC with Gurobi first, followed by AC with Ipopt. In the MILP the non-convex constraints are relaxed. In order to get as close to an LP as possible, McCormick envelopes were used for the bilinear terms in the quadratic constraints for active power. This comes from the branch flow limit and the resulting SOCP relaxation was added to the UC model. It led to about a 75\% solve time reduction and was easier to solve than the standard convex relaxation. No relaxation was used for ACOPF.

Further techniques applied were slack variables for sparse constraints, reduction of lengthy linear expressions, reformulated downtime-dependent start-up costs, and min/max constraint reformulation. The control horizon was broken down to a single time step. A single upper/lower bound pair was sought that satisfied all time steps. It limited accuracy and was only partly used for datasets with more than 2000 buses. For the smaller networks the full formulation was applied.

The team worked entirely in Python. In the UC part Pyomo was used and for ACOPF the direct interface CYIPOPT. Algorithms were not run in parallel but all available threads were used for Gurobi. Ipopt was run in serial mode.

Many details of the approach used were finalized after intense testing. This was done with respect to inclusion of reserves and of various simplifications. UC decisions that were suboptimal or infeasible were not updated. Luckily in most cases feasibility was preserved. However, a feasible solution will always exist due to the addition of slack variables.}

\subsubsection{Performer D}
\label{sec:d}

{The team, different from most others, implemented in C. It uses unit status results from a fast UC solver, schedules from multi-period ACOPF, then generates constraints if there are infeasibilities as from base case or contingency limits. The fast UC solver is Lagrangian based. The nonlinear solver is an own IPM code which utilizes Pardiso for linear algebra.

Unimportant constraints such as the subset of units providing reserves were eliminated. Upper and lower bounds were determined to decouple the multi-period ACOPF. Various techniques were used in Pardiso to speed up the solution of a single-period ACOPF such as choice of the initial point, sparse matrix techniques, step size choice, scaling. Post processing was required to improve or verify the results. A fast-decoupled AC model, factorizations, and reserve optimization were utilized.

If UC decisions needed to be adjusted, a linear sensitivity factor is generated from the violated constraint and sent back to the UC model. Typically only a few iteration were needed to substantially reduce the power mismatch.

HPC was leveraged to calculate shift factors for quick power flow solutions. This way thousands of contingencies could be evaluated very quickly.
In order to reduce congestion and line losses, line switching was implemented through selection of line candidates and the use of a linearized AC network model. Very little benefit was observed for the reduction of line losses but some benefit could be expected  by improving power flow convergence.
}

\subsubsection{Performer E}
\label{sec:e}

{This team did not utilize any discrete solver. Solution consisted of two NLP solves. Python was used to program and then AMPL as modeling language. Connection between the two was very costly timewise but KNITRO provides a good IPM implementation for which AMPL's AD efficiently delivers the derivatives. Substitution of penalty variables was important for reducing the size of the problem.

The main steps of the two-step algorithm are to first fix some integer variables, to relax the remaining integer variables, and then to solve the resulting NLP to low precision. The remaining integer variables were rounded and fixed. The rounding was no simple rounding. Feasibility was encouraged and a conservative approach taken.

Next, the NLP was resolved to high precision warm-starting from the previous solution. If no feasibility was obtained a repair was attempted.

In addition, it was tried to optimize over all time periods with some simplification for the largest networks. Overall, many parts of the algorithm could be executed rigorously. The smaller networks could be solved with KNITRO as MINLP. This is due to the fact that KNITRO has effective heuristics for the discrete variables resulting in good but in general suboptimal solutions in relatively short time. On the other hand, very large continuous NLPs with about 10m variables can be solved with KNITRO as well. 

For larger networks linear approximations were used to balance equations and limit constraints. DC power flow was not used. For the largest cases several strategies were utilized such as fixing all lines to "on", all producer/consumer binaries to priors and the shunts as well. Some continuous variables were treated as continuous across all time periods such as pon, pmax, pmin. Some reserve variables could be fixed to zero.

The infeasible solution repair hardly ever needed to be invoked. Not much attention was paid to contingencies. Essential strategy was to ensure network connectedness. Many things could be improved with additional time.}

\subsubsection{Performer F}
\label{sec:f}
  
{The key components of the algorithm are: for UC a linear approximation of the active flow plus a full horizon and for ACOPF to allow some deviation from supply/demand at a penalty. The algorithm proceeds in up to 64 "ranks" in parallel where the ranks are: Solve the UC, then alternatingly polish the UC solution, solve each ACOPF and polish the UC solution using ACOPF, predict true objective, write to file if improving. This writeout is considered a particular advantage because at any time the currently best solution is available.

There are two main reasons the results were not as good as expected, a technical one being an inexplicable memory issue that only occurred on the competition machine and the fact that transmission switching did not work well and was taken out. Very positive was the fact that all available time could be used for improvements due to the above writeout.

The team decided to modify the up and down time constraints to a convex hull approximation. Use of the ADMM algorithm with added line coverage constraints. Gurobi's lazy constraint handler is successful in dealing with the few violated voltage constraints. Smart linking of UC and ACOPF components through ADMM. 80\% of UC solve time spent on LP. Manual preprocessing of LP needed to be augmented by various heuristics for fixing variables to integers. If necessary use "unfixing" strategies.

In the ACOPF subproblems IPOPT was used with MA97 in parallel for smaller and serial for larger networks. Warm start was done including the use of appropriate options for IPOPT. UC constraints with bounded P/Q and rounded shunts by relaxing and then rounding to the nearest integer. No spatial decomposition.

For a fast contingency analysis sensitivity factors were computed using rank 1 updates of the PTDF matrix. Still it is difficult to compute the cost of all solutions. DC-based heuristic was applied to switch overloaded lines sequentially. In the future a warm-start functionality added to the LP solver promises substantial performance improvements.}

\section{Results of the GO Competition Challenge 3}
\label{sec:results}

In this section we give an overview of the competition results. In Section \ref{sec:sol_eval_scor} we explain the solution evaluation and scoring procedures, in Section \ref{sec:data} we briefly describe the test problem instances, in Section \ref{sec:general_results} we show competition-wide results, and then in Section \ref{sec:specific_results} we present results on individual solvers.

\subsection{Solution Evaluation and Scoring Procedures}
\label{sec:sol_eval_scor}

{For the sake of guaranteeing the feasibility of the problem, we employed the soft constraints rather than hard constraints for power balance, reserve balance, line limits, and post-contingency line limits. The advantage of the linear formula for the soft constraints would be its concise form that does not increase the computational complexity. However, it allows small degree of violation regardless the penalty factors. Given that the problem already include highly nonlinear and nonconvex formula, other forms might be more relevant to enforce the constraints such as logarithmic function. Nonetheless, we tested solutions identified by the performers, and announced solutions only if they passed our internal checkers. The internal checkers are designed to carefully examine if any constraints are violated.}

The solver codes implemented by the competition entrants were assessed by running them on a selection of problem instances on a common hardware platform. Each solver run took place on a dedicated node of the PNNL Deception high performance computing cluster with a fixed time limit. At the time limit, the solver was terminated if it was still running, and if a solution file was present in the working directory, then that solution file was evaluated. If no solution file was present, or if the solution file was incorrectly formatted, then the solver run was awarded a score of $0$. Solution evaluation was performed with a Python code (https://github.com/GOCompetition/C3DataUtilities) written by the GO Competition support team. The solution evaluation procedure determined whether the solution was feasible and computed the objective value. If the solution was deemed infeasible, then the solver run was awarded a score of $0$. If the solution was feasible, then the solver run was awarded a score equal to the higher of the computed objective value and $0$.

The scores of all runs of all the solvers submitted by the competition entrants on all the problem instances were used to determine a total score of each entrant, and these total scores were used to rank then entrants. The rankings were used to award prizes to the entrants. In addition to the determination of feasibility and computed objective, the evaluation code produced output on each solver run showing the largest constraint violation of each type, giving the breakdown of the objective into its various components, and recording data on the use of computational resources. This section describes the analysis performed by the GO Competition support team on the different solvers using this evaluation output.

The procedure for invoking solvers, evaluating, scoring, ranking, and awarding prizes was followed over four events. In each event, a different selection of problem instances was used for each of the three applications, and each application was evaluated as a separate division of the competition with a different time limit appropriate to the application with separate rankings and award in each division. The real time application was division 1, the day ahead application was division 2, and the week ahead application was division 3. In these three divisions, the total score of each entrant was simply the sum of their scores on all of the corresponding problem instances. Each competition event had a different set of problem instances and new awards. This section of the paper focuses on just the final event (event 4), where the largest prizes were awarded, and the competition entrants had the most time to perfect their solvers.

For the analysis of solutions, we first introduce some solver concepts. In addition to the solvers submitted by the competition entrants, the support team implemented a benchmark solver (\cite{osti_2202592}). The benchmark solver was invoked, and its solution evaluated, alongside the solvers submitted by entrants. For each problem instance, after all solvers were run and the solutions evaluated, the maximum objective value observed among all the solutions deemed feasible was called the ensemble objective, and the solution attaining that objective value was called the ensemble solution. Notionally, a solver could be implemented by running all of the submitted solvers and the benchmark solver simultaneously and taking the best resulting solution. Such a solver would then be called the ensemble solver.
This ensemble approach is comparable to the approaches of many commercial optimization solvers, for example running a barrier solver at the same time as a simplex algorithm and terminating with whichever solution is available first.

We also introduce some concepts of feasibility and objective components. Some of the constraints defined in the problem formulation are critical to even making sense of a solution. For example, the solution would be meaningless if variables that are constrained to be integers were assigned non-integer values or if real power dispatch variables were assigned values outside of the domain of definition of the cost function. These constraints and others that are relatively easy to ensure by simple projection methods are treated as hard constraints. Any solution violating these hard constraints is deemed infeasible. A numerical tolerance is used in some cases to ensure that a determination of infeasibility is not made unfairly.

A solution satisfying all these hard constraints but possibly violating the remaining soft constraints is called feasible. The soft constraints include constraints, such as power balance, that would typically need to be satisfied in a real world application, and we recognize the somewhat artificial concept of feasibility used here by describing such solutions - which may violate the soft constraints - as evaluation feasible.

The soft constraints include some constraints, such as power balance, that correspond to physical laws, and some constraints, such as line flow limits, that correspond to practical engineering limits. We therefore describe a solution satisfying all the hard constraints as well as the physical laws as physically feasible. And we describe a solution satisfying the hard constraints, and the physical law constraints, and the engineering limits as engineering feasible. A fully successful solution is one that is engineering feasible, but we have set up our evaluation procedure so that it is possible to at least evaluate physically feasible and evaluation feasible solutions as well and assess their merit alongside physically feasible solutions. We use realistic penalty coefficients to penalize violations of the soft constraints so that typically an engineering feasible solution will score better than a physically feasible solution that is not engineering feasible, and a physically feasible solution will score better than an evaluation feasible solution that is not physically feasible.

In order to understand the difference in quality between two feasible solutions, a natural approach is to observe the difference in their objective values. However, it is typically not clear how much of a difference is significant. Therefore one often considers the relative difference in objective values, i.e. as a percent of the objective value of one or the other solution or of a reference objective value. However, we observed that frequently the value of the load served is such a huge dollar amount as to dominate reasonably significant differences in solutions. Therefore, another reasonable approach is to view objective differences relative to only the generation cost.

\subsection{Note on Construction of Data and Test Problem Instances}
\label{sec:data}

Hundreds of individual problem instances were constructed to test the solvers on.
Generally this was done by constructing a few power system networks and then constructing a number scenarios for each network in each the three applications - Real Time, Day Ahead, and Week Ahead. The networks had a range of sizes, including different numbers of buses, transmission lines, transformers, phase shifters, shunts, loads, and other components. Scenarios within a network differed by conditions of weather that translated into load profiles and wind and solar availability over the model time horizon. The number of scenarios and the rough size of each network in each application is given in Table \ref{tab:net_scen}

\begin{table}[h!]
\centering
\begin{tabular}{rrrrr}
\hline
\multirow{2}{6em}{Network Size (Buses)} & \multicolumn{4}{c}{Application} \\
\cline{2-5}
 & Real Time & Day Ahead & Week Ahead & All \\
 \hline
73 & 24 & 56 & 24 & 104 \\
617 & 54 & 24 & 24 & 102 \\
1576 & 24 & 0 & 24 & 48 \\
2000 & 18 & 18 & 3 & 39 \\
4224 & 28 & 24 & 24 & 76 \\
6049 & 36 & 18 & 24 & 78 \\
6708 & 13 & 9 & 21 & 43 \\
6717 & 36 & 18 & 3 & 57 \\
8316 & 48 & 24 & 44 & 116 \\
23643 & 2 & 2 & 2 & 6 \\
\hline
All Networks & 283 & 193 & 193 & 669 \\
\hline
\end{tabular}
\caption{Networks and scenarios used in the final event.}
\label{tab:net_scen}
\end{table}

A number of the problem instances in the Division 3 week ahead application were designed to model extreme weather conditions through the load and variable generation parameters. We did not specifically study the performance of the submitted solvers on these instances in detail. However we do note that these instances were solved successfully by the submitted solvers and the benchmark solver, along with the other instances, demonstrating the ability of the solvers developed in this competition to address week ahead grid planning for extreme weather.

\subsection{Competition Results Over All Solvers}
\label{sec:general_results}

\begin{figure}[ht]
\centering
\includegraphics[width=\textwidth]{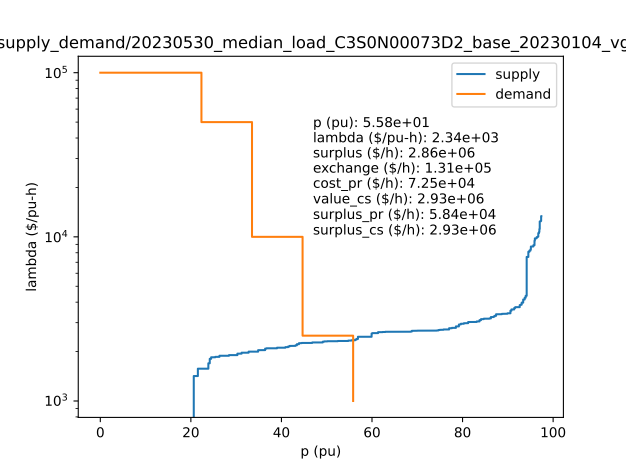}
\caption{System-wide aggregate real power supply and demand curves for a single time interval in a 73-bus problem instance.}
\label{fig:supply_demand}
\end{figure}

It is also not always clear how to obtain a reference objective value or a reference value of the generation cost. One could use values obtained from the ensemble solution or the benchmark solution, but it is not guaranteed a priori that these solutions will be particularly high quality. Therefore, we also developed a crude approximation of the optimal objective value by considering only the convex features of the value of load bids and the cost of generation offers. Features of local power balance, ramping, unit commitment, losses, reactive power, and voltage are ignored. This results in a single system-wide market for real power in each time interval, and the market clearing equilibrium can be easily computed to obtain a market surplus value that can be used as a reference objective value. The cost of generation in this market clearing equilibrium can also be used as a reference generation cost value. The market clearing equilibrium computed in this way behaves as a relaxation of full problem under certain assumptions including no negative price offers. In general, the problem instances created by the GO Competition do not satisfy these assumptions, so the market clearing equilibrium is not guaranteed to be a relaxation. The system-wide aggregate supply and demand curves are plotted in Figure \ref{fig:supply_demand} for a single time interval in one of the competition problem instances. The equilibrium real power dispatch and price are represented by the intersection of the two curves, and the area to the right of the intersection and between the two curves represents the equilibrium market surplus.

\begin{figure}[ht]
\centering
\includegraphics[width=\textwidth]{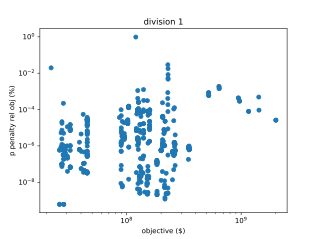}
\caption{Real power imbalance penalty as a percent of total objective vs. total objective in division 1 ensemble solutions.}
\label{fig:p_pen_vs_obj_div_1_ensemble}
\end{figure}

\begin{figure}[ht]
\centering
\includegraphics[width=\textwidth]{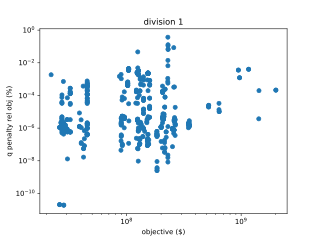}
\caption{Reactive power imbalance penalty as a percent of total objective vs. total objective in division 1 ensemble solutions.}
\label{fig:q_pen_vs_obj_div_1_ensemble}
\end{figure}

\begin{figure}[ht]
\centering
\includegraphics[width=\textwidth]{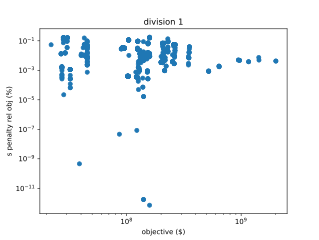}
\caption{Base case line overload penalty as a percent of total objective vs. total objective in division 1 ensemble solutions.}
\label{fig:s_pen_vs_obj_div_1_ensemble}
\end{figure}

Now, to analyze the results of the solvers developed by competition entrants, we begin with feasibility. First we examine the significance of soft constraint penalties as a way of understanding the degree to which the solvers satisfied these constraints. Figure \ref{fig:p_pen_vs_obj_div_1_ensemble} shows the value of the real power imbalance penalties, relative to the total objective, plotted against the objective value, over all problem instances in Division 1. Each point plotted represents the ensemble solution on a single problem instance. We see that the real power imbalance penalties were not significant, generally less than 1\% and typically much less. Penalties on reactive power imbalance and on line overload in the base case, plotted in Figures \ref{fig:q_pen_vs_obj_div_1_ensemble} and \ref{fig:s_pen_vs_obj_div_1_ensemble} were similarly insignificant. Therefore, the ensemble solutions were largely feasible with respect to constraints representing physical laws and with respect to the line limit engineering constraints. Similar observations hold for Divisions 2 and 3.

\begin{figure}[ht]
\centering
\includegraphics[width=\textwidth]{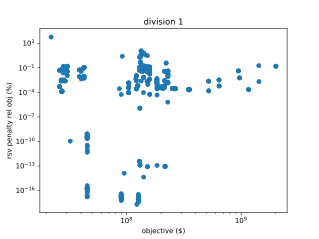}
\caption{Reserve imbalance penalty as a percent of total objective vs. total objective in division 1 ensemble solutions.}
\label{fig:rsv_pen_vs_obj_div_1_ensemble}
\end{figure}

Reserve imbalance penalties were somewhat more significant, often around 1\% to 10\%, as shown in Figure \ref{fig:rsv_pen_vs_obj_div_1_ensemble}. One way to view the reserve imbalance penalties is that reserve requirements are often modeled with a reserve demand curve, so that the requirement is merely one point on the curve, and these reserve imbalances correspond to situations of reserve scarcity.

\begin{figure}[ht]
\centering
\includegraphics[width=\textwidth]{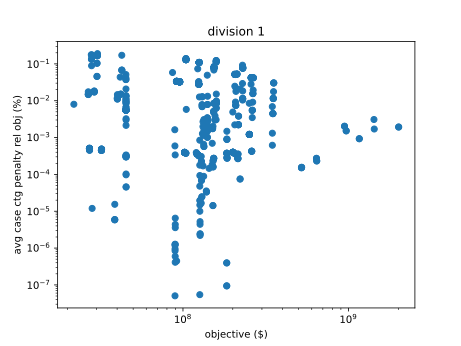}
\caption{Average case post-contingency line overload penalty as a percent of total objective vs. total objective in division 1 ensemble solutions.}
\label{fig:avg_post_ctg_pen_vs_obj_div_1_ensemble}
\end{figure}

\begin{figure}[ht]
\centering
\includegraphics[width=\textwidth]{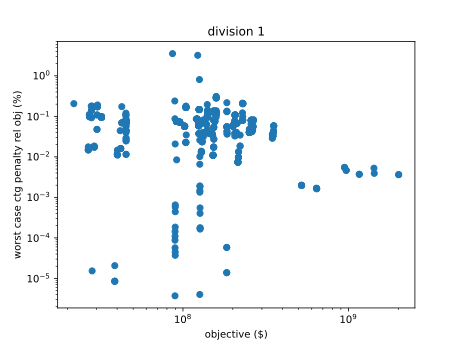}
\caption{Worst case post-contingency line overload penalty as a percent of total objective vs. total objective in division 1 ensemble solutions.}
\label{fig:worst_post_ctg_pen_vs_obj_div_1_ensemble}
\end{figure}

Post-contingency line overload penalties were similarly insignificant, as shown in Figures \ref{fig:avg_post_ctg_pen_vs_obj_div_1_ensemble} (average case) and \ref{fig:worst_post_ctg_pen_vs_obj_div_1_ensemble} (worst case). As expected, the worst case penalties, on the order of up to 1 \% in a few instances, were more significant than the average case penalties, but still relatively small. The worst case penalty term was included in order to ensure that solvers face a nontrivial incentive to respect the post-contingency line constraints.

\begin{figure}[ht]
\centering
\includegraphics[width=\textwidth]{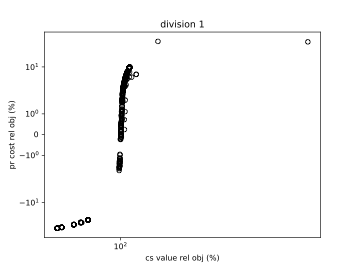}
\caption{Generator production cost as a percent of total objective vs. consumer load value as a percent of total objective in division 1 ensemble solutions.}
\label{fig:pr_cost_vs_cs_value_div_1_ensemble}
\end{figure}

\begin{figure}[ht]
\centering
\includegraphics[width=\textwidth]{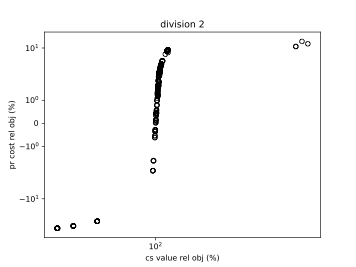}
\caption{Generator production cost as a percent of total objective vs. consumer load value as a percent of total objective in division 2 ensemble solutions.}
\label{fig:pr_cost_vs_cs_value_div_2_ensemble}
\end{figure}

\begin{figure}[ht]
\centering
\includegraphics[width=\textwidth]{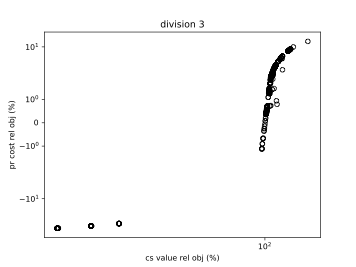}
\caption{Generator production cost as a percent of total objective vs. consumer load value as a percent of total objective in division 3 ensemble solutions.}
\label{fig:pr_cost_vs_cs_value_div_3_ensemble}
\end{figure}

We can perform similar analysis to determine the most significant components of the ensemble objectives. In Figure \ref{fig:pr_cost_vs_cs_value_div_1_ensemble}, the generator production cost as a percent of the total objective value is plotted against the consumer load value as a percent of the total objective for Division 1. Each plotted point corresponds to the ensemble solver on a single problem instance. This plot shows that the consumer value is typically between 90\% and 110\% of the total objective, the generator cost is typically between -10\% and 10\% of the objective. Furthermore, we see that gains in load value are realized in a one-for-one fashion at the expense of increased generator cost. Figures \ref{fig:pr_cost_vs_cs_value_div_2_ensemble} and \ref{fig:pr_cost_vs_cs_value_div_3_ensemble} show the corresponding plots for Divisions 2 and 3 and we make similar observations there.

\begin{figure}[ht]
\centering
\includegraphics[width=\textwidth]{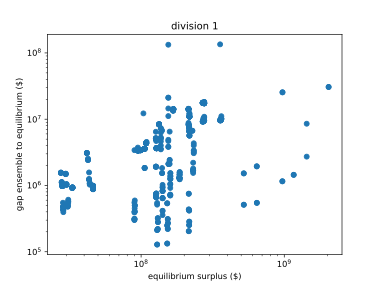}
\caption{Absolute gap of objective to equilibrium market surplus vs. equilibrium market surplus in division 1 ensemble solutions.}
\label{fig:absolute_gap_vs_equil_div_1_ensemble}
\end{figure}

\begin{figure}[ht]
\centering
\includegraphics[width=\textwidth]{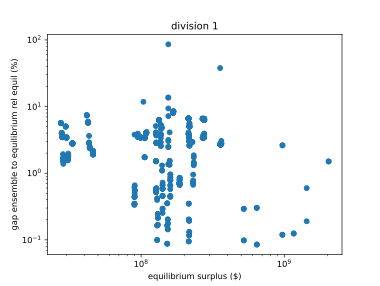}
\caption{Gap of objective to equilibrium market surplus relative to equilibrium market surplus vs. equilibrium market surplus in division 1 ensemble solutions.}
\label{fig:gap_rel_equil_vs_equil_div_1_ensemble}
\end{figure}

\begin{figure}[ht]
\centering
\includegraphics[width=\textwidth]{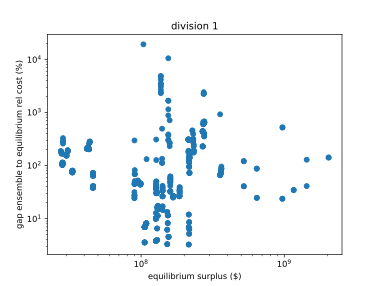}
\caption{Gap of objective to equilibrium market surplus relative to equilibrium generator production cost vs. equilibrium market surplus in division 1 ensemble solutions.}
\label{fig:gap_rel_cost_vs_equil_div_1_ensemble}
\end{figure}

The analysis of penalty values shows that the ensemble solutions were generally feasible. To assert that the ensemble solutions fully solved the competition problem we need to assess optimality. To rigorously prove that a feasible solution is optimal, we need to compare its objective value to that of a relaxation and show that the gap is small enough for practical use. We do not have a theoretically rigorous relaxation procedure, but the market clearing equilibrium is a relaxation under certain conditions and appears to behave like a relaxation. Figure \ref{fig:absolute_gap_vs_equil_div_1_ensemble} plots the absolute gap from the ensemble objective value to the equilibrium market surplus value against the equilibrium market surplus value in Division 1. In all cases the gap is positive, showing that we have not found a case where the relaxation behavior of the market equilibrium does not hold. In Figure \ref{fig:gap_rel_equil_vs_equil_div_1_ensemble}, we normalize the gap by the equilibrium market surplus value, and we see that most of the time the gap is less than 10\% and frequently less than 1\%. This supports the claim that the ensemble solutions are nearly optimal. In Figure \ref{fig:gap_rel_cost_vs_equil_div_1_ensemble}, we normalize the gap by the generator production cost in the market clearing equilibrium solution instead of by the total objective, and by this measure the gap is typically on the order of 10\% to 1000\%. With this more stringent normalization, we see that there might very well be significant room for improvement in some of the solutions, or else the relaxation bound provided by the market equilibrium solution might not be very tight, or some combination of the two. This shows that interpretation of the gap between two solutions depends rather strongly on the method of normalization.





\subsection{Competition Results Specific To Individual Solvers}
\label{sec:specific_results}

Table \ref{tab:team_performance} presents an overview of the results of
the GO Competition Challenge 3 for the individual solvers submitted by the competition entrants.

In the table, each team is shown on a single row. The first row represents a hypothetical ensemble solver that gives the best solution of all the teams on any given scenario. The second row is the ARPA-e Benchmark solver, which was not eligible for prize awards. Thereafter the teams are sorted in order of decreasing final event prize award.

Columns labeled Obj represent the total objective value (market surplus minus constraint violation penalties) achieved by the given team over all scenarios, in \$1e9. Columns labeled NB represent the number of scenarios on which the given team achieved the best (highest) objective value of all teams. 

For each of the three applications - Real Time, Day Ahead, and Week Ahead - and for all three combined, the total objective and number of best scenarios are reported. Also, the prize award in the final event and that in all events of the competition are reported.

The table shows that YongOptimization was the clear overall winner. TIM-GO, GOT-BSI-OPF, Artelys\_Columbia, and GravityX also gave strong performances. The Blackouts, Occam's Razor, and Electric-Stampede made significant contributions. And we know that LLGoMax, quasiGrad, Gatorgar, PACE, and PGWOpt put in a great deal of effort. We also point out that the ARPA-e Benchmark solver achieved a very respectable total objective score.
YongOptimization found the best solution on the greatest number of scenarios, but TIM-GO and GOT-BSI-OPF achieve a higher total objective score, presumably reflecting better performance on certain very large problem instances. And the fact that the teams other than YongOptimization achieved a better solution on about one half of the scenarios shows that many teams made unique contributions to solving this challenging problem. This high performance achieved by a diversity of approaches is precisely the kind of result that can be achieved by a competition approach to research, as it incentivizes researchers to try radically different methods.

\begin{table}[h!]
\centering
\begin{tabular}{rrrrrrrrrrr}
\hline
& \multicolumn{2}{c}{Real} & \multicolumn{2}{c}{Day} & \multicolumn{2}{c}{Week} & \multicolumn{2}{c}{} & \multicolumn{2}{c}{Award} \\
& \multicolumn{2}{c}{Time} & \multicolumn{2}{c}{Ahead} & \multicolumn{2}{c}{Ahead} & \multicolumn{2}{c}{All} & \multicolumn{2}{c}{(\$1e3)} \\
\cline{2-3} \cline{4-5} \cline{6-7} \cline{8-9} \cline{10-11}
Team & Obj & NB & Obj & NB & Obj & NB & Obj & NB & Final & All \\
\hline
Ensemble & 45.8 & 283 & 163.5 & 192 & 914.9 & 192 & 1124.4 & 667 & 0 & 0 \\
ARPA-e Benchmark & 40.7 & 0 & 156.0 & 0 & 894.7 & 0 & 1091.4 & 0 & 0 & 0 \\
YongOptimization & 44.5 & 156 & 160.1 & 78 & 898.4 & 97 & 1103.1 & 331 & 550 & 645 \\
TIM-GO & 43.8 & 39 & 162.2 & 23 & 912.9 & 27 & 1119.1 & 89 & 520 & 595 \\
GOT-BSI-OPF & 45.1 & 28 & 162.9 & 15 & 912.2 & 4 & 1120.3 & 47 & 360 & 485 \\
Artelys\_Columbia & 41.9 & 2 & 157.3 & 17 & 890.9 & 10 & 1090.2 & 29 & 320 & 390 \\
GravityX & 41.0 & 37 & 156.1 & 38 & 615.7 & 30 & 812.9 & 105 & 320 & 440 \\
The Blackouts & 16.2 & 16 & 114.0 & 18 & 460.2 & 22 & 590.6 & 56 & 200 & 200 \\
Occams razor & 42.0 & 3 & 145.4 & 1 & 859.3 & 0 & 1046.8 & 4 & 130 & 130 \\
Electric-Stampede & 33.7 & 0 & 139.3 & 0 & 789.5 & 0 & 962.6 & 0 & 0 & 115 \\
LLGoMax & 20.9 & 1 & 116.8 & 0 & 743.7 & 2 & 881.4 & 3 & 0 & 0 \\
quasiGrad & 31.4 & 0 & 155.1 & 0 & 507.6 & 0 & 694.3 & 0 & 0 & 0 \\
Gatorgar & 2.4 & 1 & 10.2 & 2 & 43.0 & 0 & 55.7 & 3 & 0 & 0 \\
PACE & 0.0 & 0 & 0.0 & 0 & 0.0 & 0 & 0.0 & 0 & 0 & 0 \\
PGWOpt & 0.0 & 0 & 0.0 & 0 & 0.0 & 0 & 0.0 & 0 & 0 & 0 \\
\hline
\end{tabular}
\caption{Performance of solvers in the final event}
\label{tab:team_performance}
\end{table}

\section{Summaries and Future Research Directions}
\label{sec:conclusion}
{GO Challenge 3 anticipates an era marked by a harmonious blend of sophisticated computational techniques and nuanced, adaptable algorithmic strategies in the future of UC and AC OPF research. This vision is the culmination of the collective insights of teams that have pioneered various approaches in these fields.

The concept of dynamic, scalable modeling is central to this vision. Rolling horizon approaches and the combination of spatial and temporal decomposition are being pioneered by teams such as Occam's Razor and GravityX. These strategies represent a shift toward models that can adapt to the ever-changing dynamics of power systems, allowing them to accommodate large-scale datasets and fluctuating conditions with greater precision and efficiency.

In addition, there is a clear trend toward combining advanced computational methods with traditional mathematical optimization. This is evident in the GravityX team's interest in combining machine learning and mathematical optimization, pointing to a future in which predictive analytics will play a critical role in improving the accuracy and efficiency of power system management. This type of integration could result in algorithms that not only solve current problems more effectively, but also anticipate future challenges, resulting in more resilient power systems.

Another critical theme is efficiency in algorithmic design. Electric Stampede and Artelys have shown the importance of tailoring computational strategies to the specific needs of different power networks, whether through sophisticated handling of complex constraints or balancing computational efficiency with solution accuracy. This indicates a future where solutions are not one-size-fits-all but are instead customized to the unique characteristics of each power system network.

Furthermore, the Yong Optimization team's exploration of High-Performance Computing (HPC) and advanced numerical methods highlights the growing importance of computational power in managing complex power systems. This suggests that greater computational resources are being used to solve power system problems more quickly and accurately.

Finally, as seen in The Blackouts' approach, the integration of parallel processing techniques and the development of new computational tools points to a future in which solving UC and ACOPF problems becomes more efficient and effective. Incorporating warm-start capabilities into LP solvers, for example, is a step toward reducing computation times and improving the overall performance of optimization algorithms.

In conclusion, the future of UC and AC OPF research is shaping up to be an interdisciplinary endeavor combining dynamic modeling, computational intelligence, algorithmic efficiency, and high-powered computing. This integrated approach is poised to more effectively address the complexities of modern power systems, paving the way for more robust, adaptable, and efficient power system management solutions.
}
\section*{Disclaimer}

 The views, information, and opinions expressed in this paper are solely those of the authors and do not necessarily represent those of the GO Competition performers, PNNL, ARPA-E, DOE, or the United States Government.
 The authors have no conflicts of interest.

%

\printbibliography

\end{document}